\documentclass[12pt]{mcmthesis}
\mcmsetup{tstyle=\color{red}\bfseries,
        tcn = 2401540, problem = D, 
        sheet = false, titleinsheet = false, keywordsinsheet = true,
        titlepage = true, abstract = true}

  \usepackage{newtxtext,newtxmath} 
  \usepackage{graphicx}
  \usepackage{tocloft}

\makeatletter
\let\c@lofdepth\relax
\let\c@lotdepth\relax
\makeatother
\usepackage{subfigure}	
\usepackage{float}	
\usepackage{enumitem}
\usepackage{url}
\usepackage{algorithm}
\usepackage{algorithmicx}
\usepackage{algpseudocode}
\usepackage{cite}
\usepackage{indentfirst}  
\usepackage{lipsum}

\title{WLPCM Approach for Great Lakes Regulation}
\author{Wenbo Huang, Xiangyi Chen, Jiaqi Leng}
\date{February 3, 2024}
\begin{document}

\begin{abstract}
\par 
Managing the water levels of the Great Lakes, the largest freshwater lakes globally, is a great challenge. Stakeholders have diverse request. This study aims to develop a model for more effective water level management, solving the challenging network flow problem in a practical sense.


For requirement 1, we exclude abnormal data and get average water levels based on data before 2019. Monthly average water levels of Great Lakes display a consistent annual pattern. With multi-year averages as a benchmark, we then aim to maximize the fitness function by controlling flow from Compensating Works and Moses-Saunders Dam, thus identifying optimal water levels. In such circumstances, we use the Simulated Annealing algorithm to find the the optimized flow.

For requirement 2, we first determine the coefficients in the basic water level model through linear function and calculate natural variables through the given data. Then, we employ the former natural Indicator to predict the latest Indicator through linear regression, using Delay Differential Equation (DDE) for dam regulation effects. Finally, we present the Water Level Predictive Control Model (WLPCM) based on Model Predictive Control (MPC). The model integrates a predictive component using a DDE-based model, considering inflow, outflow and natural variables. Due to the temporal lag of upstream changes, the model optimizes dam flow rates in six months. The model also includes feedback adjustment, ensuring a rational response, even in extreme conditions. 

For requirement 3, we tested efficacy of WLPCM based on data from 2017, utilizing Validation-DDE to simulate outcomes with altered control. Our approach effectively mitigated the unusually high water levels in the Ottawa River observed in 2017.Our scores significantly surpassed the records of 2017 in lake water level. Our model receives feedback from real world, allowing for adjustments and corrections. Employing a Sobol Approach to analyze the sensitivity of lake-water regulation, the results affirm the robustness of our model.

For requirement 4, we have made slight alterations to certain stages and evaluated the RMSE of the fitness function to assess sensitivity. The natural variations in ice-clog, snow-pack, and precipitation can be quantified in terms of a $\sigma$ change in water flow or lake level, with sensitivity assessed through the variance in grading. Our findings indicate that ice-clog has the most significant impact on the Great Lakes system among all the indicators considered.

For requirement 5, we catalogued demands of stakeholders around Lake Ontario. After analyzing the hydrodynamic models and water usage requirements, we developed a targeted Fitness Function. We refined our control strategy using WLPCM, subsequently validating the appropriateness of our control measures through DDE-based model. Our control outcomes for 2017 demonstrate that our model successfully reduced the catastrophic high water levels in Lake Ontario and the excessive flow rates in the St. Lawrence River. Compared to the Plan 2014 in use at the time, our approach effectively considered optimal long-term strategies, thereby preventing flooding.

\end{abstract}
\maketitle
\tableofcontents
\newpage
\section{Introduction}
\subsection{Background}
The Great Lakes of the United States and Canada, form the world's largest freshwater lake system, crucial for diverse purposes like fisheries, recreation, and power generation. Managing water levels is essential but challenging due to complex factors and conflicting interests. The region faces issues like floods from insufficient discharge and economic hindrance from excessive water release. Two key control mechanisms are the Compensating Works in Sault Ste. Marie and the Moses-Saunders Dam in Cornwall. Despite human control, factors like rainfall, evaporation, and ice blockages pose challenges. Effective management requires balancing natural and human influences to sustain the ecosystem and meet stakeholders' needs.


\subsection{Restatement of the Problem}
Through in-depth analysis and research on the background of the problem, the restate of the problem can be expressed as follows:
\begin{itemize}
\setlength{\topsep}{1ex} 
\setlength{\itemsep}{0.5ex} 
    \item Determine the ideal water levels for the Great Lakes at any time during the year, considering the desires of various stakeholders, each with potentially different costs and benefits.
    \item Develop an algorithm to maintain the optimal water levels for the Great Lakes based on inflow and outflow data.
    \item Evaluate the sensitivity of the control algorithm to the outflow rates of the two controlling dams. With data from 2017, assess whether the new control would result in satisfaction for stakeholders or achieve better outcomes than the recorded water levels.
    \item Examine the algorithm in respond to changes in environmental conditions (e.g., precipitation, winter snowpack, ice jams)
    \item Focus the comprehensive analysis specifically on the stakeholders and influencing factors related to Lake Ontario.
\end{itemize}

\section{Modelling and Natural Scenario}

\subsection{Assumptions}
\begin{enumerate}
\item The change in water volume of each lake is determined by inflow, outflow, evaporation, and rainfall. We neglect the extra loss during transmission, indicating that human activities have no influence on the total volume of lake water level.
\item The river flow and lake water level possess linear interrelationship, which is validated through the linear fitting process.
\end{enumerate}

\subsection{Notations}
\begin{center}
\begin{tabular}{cc}
 \toprule[1.5pt]
 \makebox[0.3\textwidth][c]{Notations}	&  \makebox[0.5\textwidth][c]{Indication} \\
 \midrule[1pt]
 $ H,F $ & Current Water Level/Flow \\
 $ H^*,F^* $  &  Average Water Level/FLow over 12 months \\
 $ \hat{H},\hat{F}$ & Original Water Level/Flow\\
 $ \hat{H^*},\hat{F^*} $ & Original Average Water Level/Flow over 12 months\\
 $ \sigma $ & Standard Deviation\\
 $ R $ & Restraints on different lakes\\
 $ \alpha_{ij}, \beta_{ij}, \gamma_{ij}, \delta_{ij}, \eta_{ij}, \Delta_k $ & Coefficients in Differential Equations\\
 $ S $ & Sensitivity\\
 $\Delta$ & Natural Indicator\\
\bottomrule[1.5pt]
\end{tabular}
\end{center}

\subsection{Modelling of water level}
\begin{figure}[H]
\centering
\includegraphics[width=0.9\textwidth]{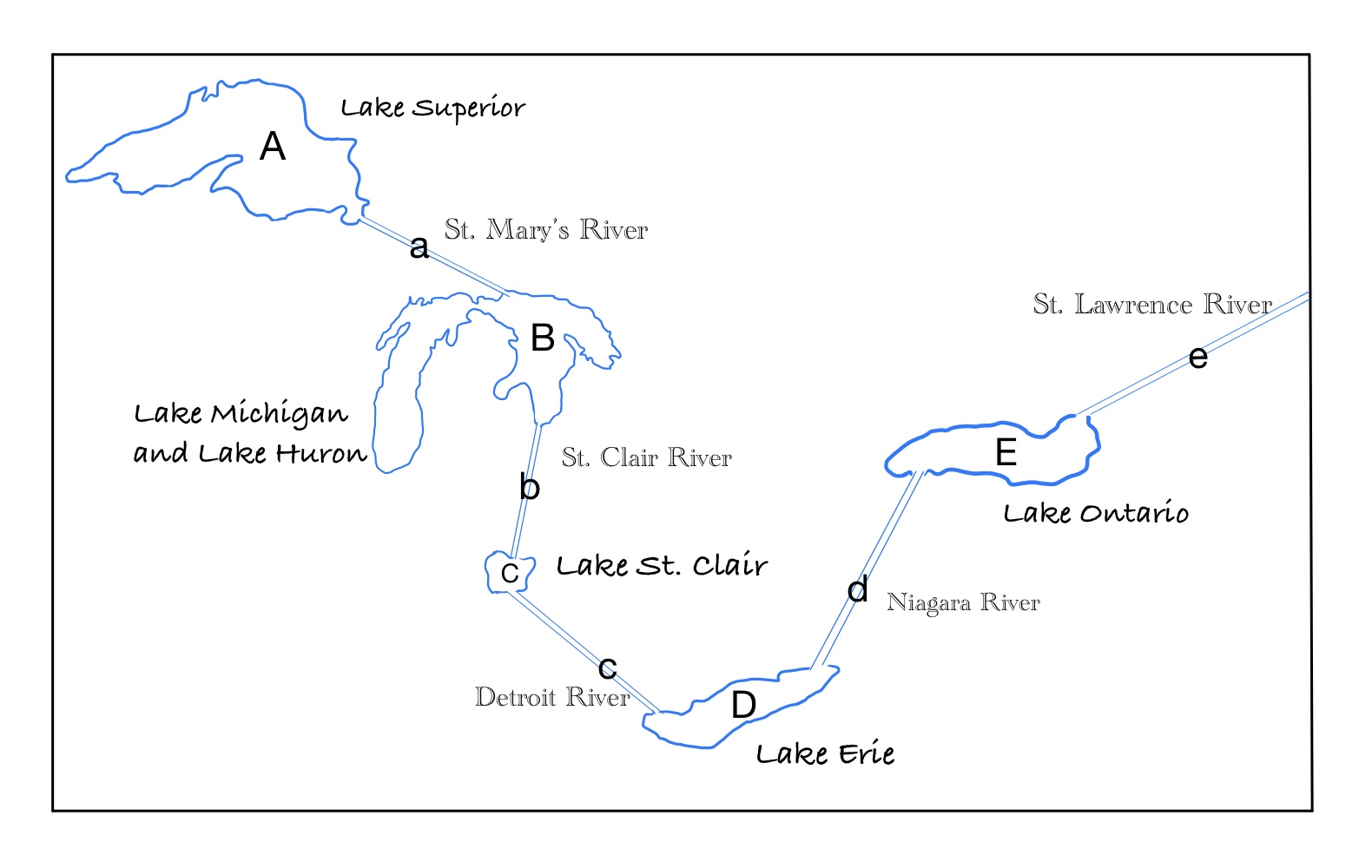}
\caption{Simplified Diagram of the Great Lakes and major Rivers}
\end{figure}
The connection between the Great lakes and major rivers can be simplified as the caption above. The water flows from west to east, and the altitude drop between different lakes is the main reason for the disparity of water level and river flow. The Montreal City locates by the St. Lawrence River and represents the outflow of the Great Lakes.

 To derive an overall relationship between the levels of the Great Lakes and the flow of water under natural conditions, we use data on water levels of the Great Lakes and the flow of rivers for past months to develop a mathematical model. The change of the water level is determined by natural factors and the water flow inwards and outwards. The set of lake functions can be interpreted as

$$
\left\{
\begin{aligned}
 &\frac{\mathrm{d} A}{\mathrm{d} t} =-\alpha_{11}a(t-1) +\Delta _{A}(t)\\
 &\frac{\mathrm{d} B}{\mathrm{d} t} =\beta _{11}a(t-1) -\beta _{12}b(t-1)+\Delta _{B}(t)\\
 &\frac{\mathrm{d} C}{\mathrm{d} t} =\gamma _{11}b(t-1) -\gamma _{12}c(t-1)+\Delta _{C}(t)\\
 &\frac{\mathrm{d} D}{\mathrm{d} t} =\delta_{11}c(t-1) -\delta _{12}d(t-1)+\Delta _{D}(t)\\
 &\frac{\mathrm{d} E}{\mathrm{d} t} =\eta _{11}d(t-1) -\eta _{12}e(t-1)+\Delta _{E}(t)\\
\end{aligned}
\right.
$$

The $\Delta$ here represents the natural factors like evaporation and osmosis that influence the levels of the Great Lakes. $A$,$B$, $C$, $D$ and $E$ respectively represent Lake Superior, Lake Michigan and Lake Huron, Lake St. Clair, Lake Erie and Lake Ontario. $a$, $b$, $c$, $d$ and $e$ respectively represent St. Mary's River, St. Clair River, Detroit River, Niagara River and Ottawa River.

The flow of river $a$, $b$, $c$, $d$ is subject to the natural factors like water flow drop and local climate, and is also influenced by the water level of the lake that the river derived from, which can be linearly fitted. River $e$ is closely related to Montreal city's ecology, and the water flow is artificially adjusted.

$$
\left\{
\begin{aligned}
 &a(t)=\alpha A(t-1)+\Delta _{a}\\
 &b(t)=\beta B(t-1)+\Delta _{b}\\
 &c(t)=\gamma C(t-1)+\Delta _{c}\\
 &d(t)=\delta D(t-1)+\Delta _{d}\\
\end{aligned}
\right.
$$

In later analysis, we can change the flow of river a and e to get a further relationship between changes in lake level and artificial adjustments.
\subsection{Fitness function}
In the basic simulation situation, we assume that the stakeholder only makes requests for lake surface water levels, and each lake surface is subject to at most one constraint on water surface level and water surface fluctuation degree respectively. Constraints can be divided into the following two categories:

\subsubsection{Water surface height}
The requirements for water surface height encompass three primary conditions: high water level, low water level and regular water level.

In the hypothesis, maritime transport entities exhibit a preference for maximally elevated water levels, whereas the local populace advocates for minimal water levels. Individuals who have no specific demand for this are regarded as regular water level demands.

Here, average water level over twelve months is employed to quantify the water level: A decline in the average water level relative to natural baselines signifies water level decrement; whereas an increment indicates a rise. 
According to the compliance of the average water level change with the predefined constraint conditions, the water surface height score of the lake is obtained.

The Grading of water level can be defined as
$$
G_{L}=\left\{
\begin{aligned}
&2+9(H-H^*);   & High-level-required \\
&4-18|H-H^*|; &  Medium-level-required \\
&2-9(H-H^*); &  Low-level-required
\end{aligned}
\right.
$$

and $H^*$ indicates the average water level in this situation:
$$
H^*=\frac{\sum\limits_{i=1}^{12}H_i}{12}
$$

The original $H_i$ would be acquired in the follow-up procedure of data processing. The more closely the water level aligns with the specified requirements, the higher the grade presents.

\subsubsection{Water fluctuation degree} 
The requirements for the degree encompass large fluctuations, small fluctuations and regular fluctuations.In the hypothesis, sailing company owners need a stable environment to assure the smooth running of their industry, so they prefer stable water level. However, environmentalists prefer dynamic of lakes so that the environment can remain its diversity, and individuals with no specific demand for this are regarded as a regular fluctuating demand.

In this context, water level assessment is conducted through the analysis of the twelve-month standard deviation in water levels. An augmentation in the standard deviation relative to natural states is interpreted as a significant fluctuation, whereas a diminution in standard deviation signifies a minor fluctuation.
The Grading of water fluctuation can be defined as 
$$
G_{F}=\left\{
\begin{aligned}
&2+12(\sigma - \hat{\sigma})  ; & High-fluctuation-required \\
&4-24|\sigma - \hat{\sigma}|  ; & Medium-fluctuation-required \\
&2-12(\sigma - \hat{\sigma})  ; & Low-fluctuation-required
\end{aligned}
\right.
$$
$\sigma$ and $\hat\sigma$ stands the newest standard deviation and original standard deviation:
$$
\sigma = \sqrt{\frac{\sum\limits_{i=1}^{12}(H-H^*)^2}{12}}
;\quad
\hat\sigma = \sqrt{\frac{\sum\limits_{i=1}^{12}(\hat{H}-\hat{H^*})^2}{12}}
$$
Subsequently, the cumulative scores of $G_L$ and $G_F$ over five distinct lake according to the lake level stability and fluctuation are generated. We additionally set the limit that $G_L$ and $G_F$ is no less than $0$ and no more than $4$. Thus, The optimal scenario is identified as the condition as water flow yields the maximum score value with range $0-8$, corresponding to the most favorable water level adjustments.
$$
G(a,e,R)=\sum^{Lake_i}_{1\leq i \leq 5}(G_{L}+G_{F}) 
$$
In the context of the problem, the level and fluctuation of water can be slightly controlled through the Compensating Works and Moses-Saunders Dam. The change of $G$ is attributed to the alteration of the flow of river $a$ and $e$, so the model can be interpreted as finding the maximum of $G(a,e,R)$ and the corresponding water flow of $a,e$ with fixed $R$.


\section{Requirement 1 - Simulated Annealing Approach to Optimize Lake Water Level}
\subsection{Data processing}

In the beforehand data analysis, major data like lake water level and water flow and excludes abnormal data to ensure the stability of data are sifted through. The implement of Plan 2014 resulted in greater fluctuation of lake level, resulting in an error larger than $3\sigma$ over previous water level data in 2020 and 2022, thus we temporarily consider the major assumption of Ideal water level on the basis of data before 2019.

\subsection{Natural Constraints}
\begin{figure}[H]
  \centering
    \subfigure[Lake Superior]{\includegraphics[width=0.35\textwidth]{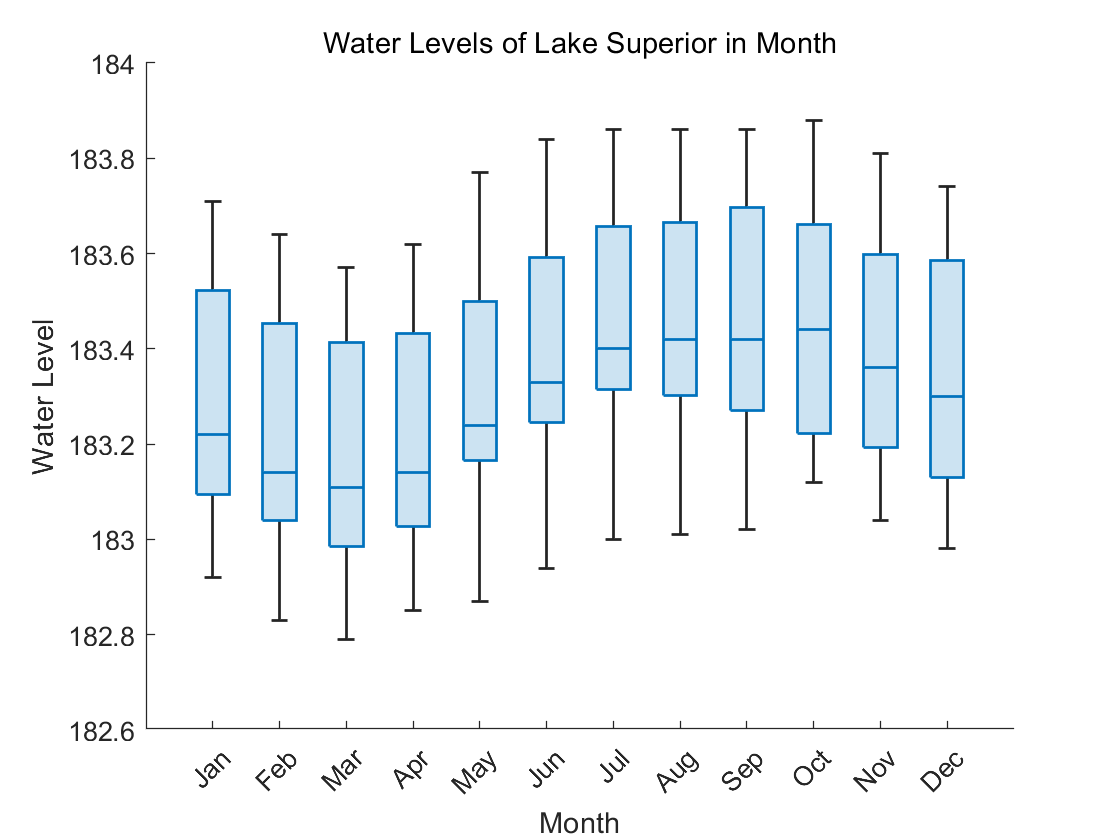}} 
    \subfigure[Lake Michigan and Lake Huron]{\includegraphics[width=0.35\textwidth]{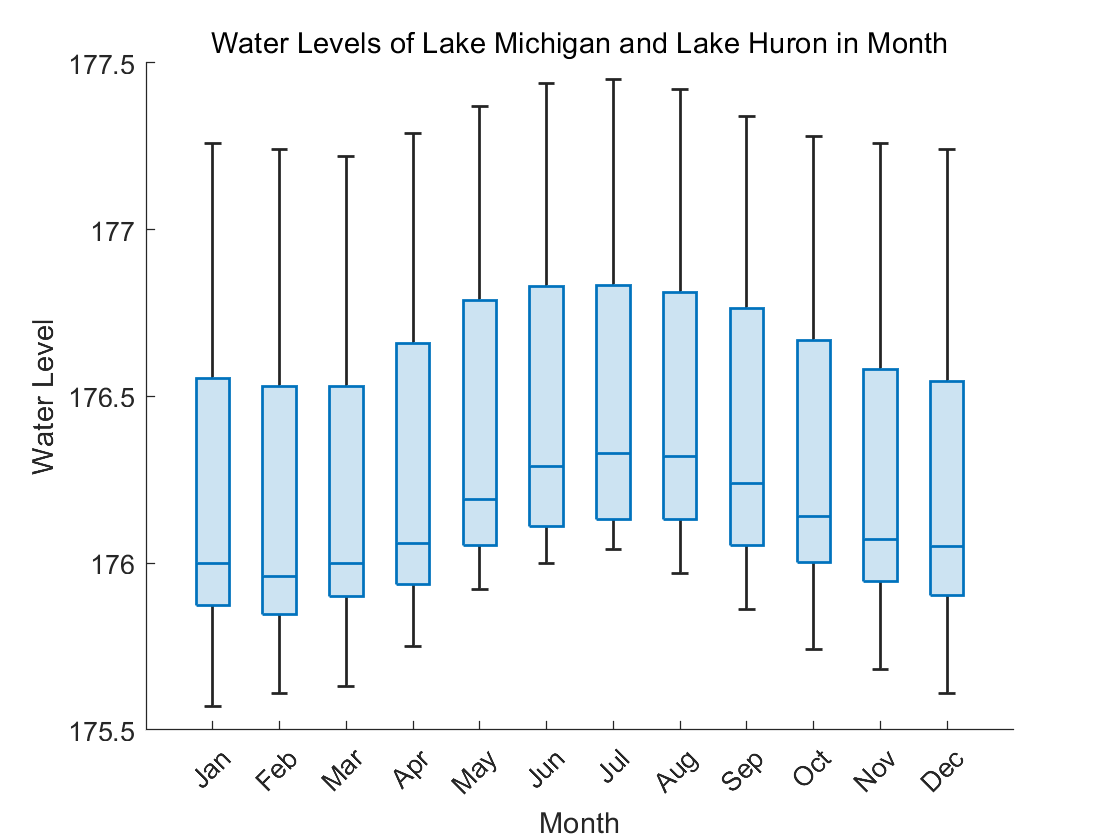}} 

    \subfigure[Lake Ontario]{\includegraphics[width=0.35\textwidth]{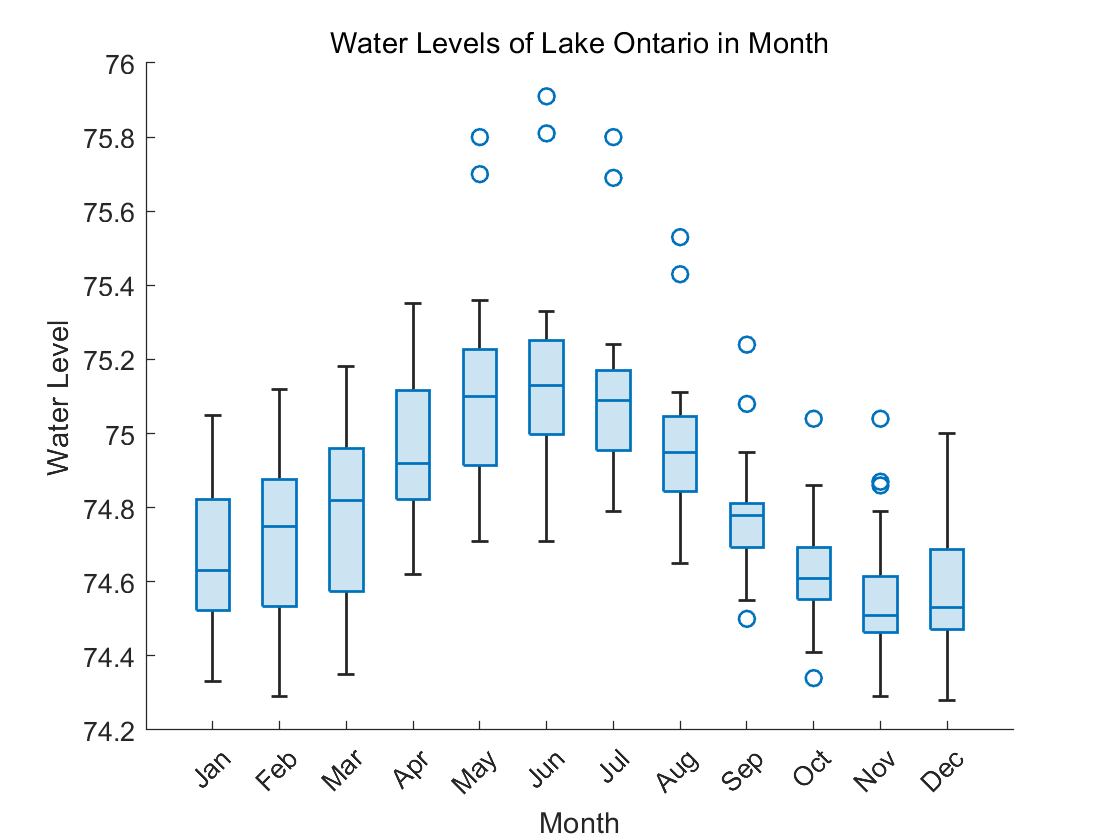}} 
    \subfigure[Lake Erie]{\includegraphics[width=0.35\textwidth]{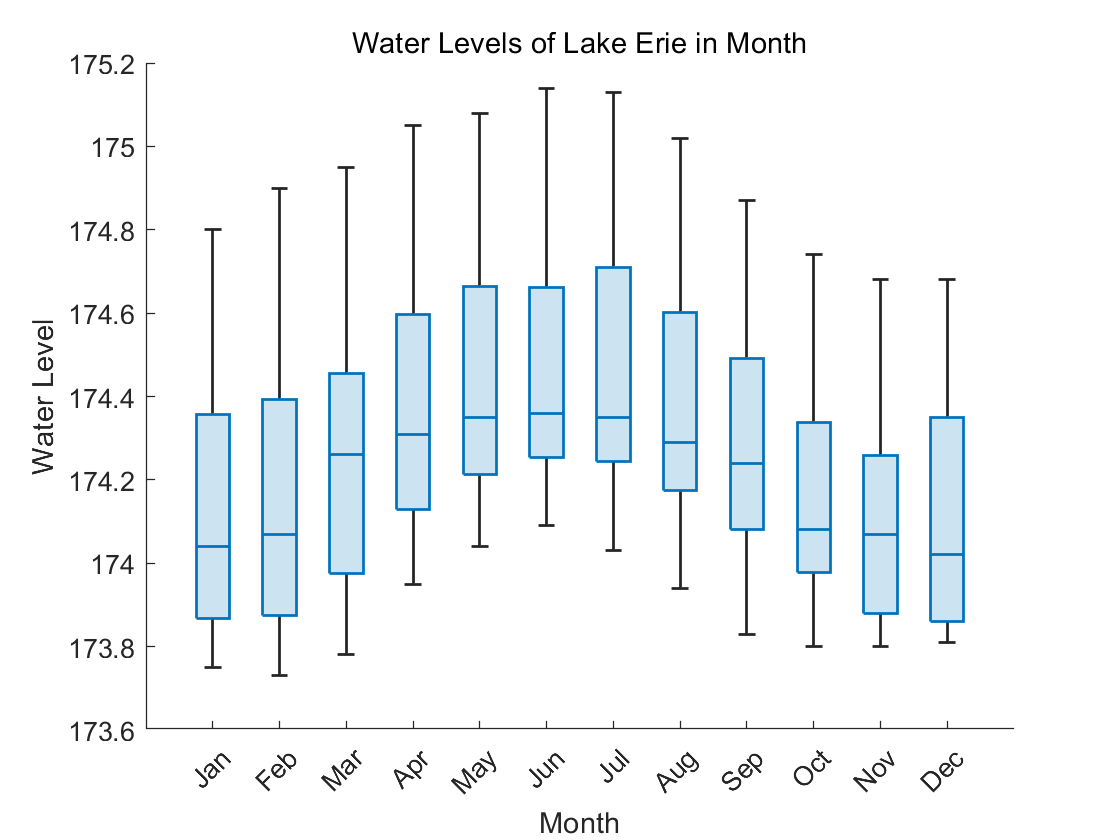}} 
    
  \caption{Boxplot of Monthly Average Lakes Water Levels}
\end{figure}

The figure above presents an analysis of the monthly average water levels of the Great Lakes. Utilizing boxplot methodology to eliminate outliers and calculate multi-year average values, we have analyzed the monthly average water levels of the Great Lakes. The flood season varies for each lake, with the occurrence of high and low water levels differing in timing, yet a consistent pattern emerges annually. This provides crucial reference points for our water management strategies for the lakes.

Given the natural conditions within the watershed, it is not always possible to meet the demands of all stakeholders. We use the multi-year average monthly data as a benchmark and consider controlling the flow from two major dams, Compensating Works and Moses-Saunders Dam, to maximize the satisfaction function. By doing so, we can identify feasible optimal water levels for the Great Lakes.

\subsection{Optimization with Simulated Annealing Algorithm}
To find the optimized flow, we employ the Simulated Annealing algorithm for solving this problem. The Simulated Annealing algorithm is a probabilistic technique for approximating the global optimum of a given function. The major algorithm can be explained in the following form:

\begin{algorithm}[H]
\caption{Simulated Annealing for Maximizing Satisfaction in Great Lakes Problem}\label{alg:simulated_annealing}
\begin{algorithmic}[1] 
\Procedure{Simulated Annealing}{}
    \State Initialize temperature $T$ to a high value.
    \State Initialize $T_{\text{min}}$, the minimum threshold.
    \State Initialize cooling rate $\alpha$.
    \State Initialize $x_{\text{current}}$ and $y_{\text{current}}$ with random 12-dimensional vectors, indicating the dam control, namely $a$ and $e$ over the whole year.
    \State Calculate $G_{\text{current}} = \text{satisfaction}(x_{\text{current}}, y_{\text{current}})$. G is the fitness function defined above.
    \State Set $x_{\text{best}}, y_{\text{best}} = x_{\text{current}}, y_{\text{current}}$ and $G_{\text{best}} = G_{\text{current}}$.
    \While{$T > T_{\text{min}}$}
        \State Generate $x_{\text{new}}$ and $y_{\text{new}}$ by making small random changes to $x_{\text{current}}$ and $y_{\text{current}}$.
        \State Calculate $G_{\text{new}} = \text{satisfaction}(x_{\text{new}}, y_{\text{new}})$.
        \State Calculate $\Delta G = G_{\text{new}} - G_{\text{current}}$.
        \If{$\Delta G > 0$ or $\exp(\Delta G / T) > \text{random}(0, 1)$}
            \State $x_{\text{current}}, y_{\text{current}} = x_{\text{new}}, y_{\text{new}}$.
            \State $G_{\text{current}} = G_{\text{new}}$.
            \If{$G_{\text{new}} > G_{\text{best}}$}
                \State $x_{\text{best}}, y_{\text{best}} = x_{\text{new}}, y_{\text{new}}$.
                \State $G_{\text{best}} = G_{\text{new}}$.
            \EndIf
        \EndIf
        \State $T = T \times \alpha$.
    \EndWhile
    \State \Return $x_{\text{best}}, y_{\text{best}}, G_{\text{best}}$.
\EndProcedure
\end{algorithmic}
\end{algorithm}

It inhibits a descending possibility to switch to the worsen solutions, allowing $G(a,e,R)$ to step out of the temporary ideal outcome while gradually focusing on areas of the search space with better solutions.

The total $G(a,e,R)$ over 12 months are the indicator of various stakeholders’ desires, and Simulated-Annealing algorithm navigates the complex landscape of possible dam flow settings to find an optimal balance that maximizes overall satisfaction through iterative refinement and the probabilistic acceptance of new solutions.

\subsection{Optimization Results}



Taking Lake Superior as an example, we present our computed outcomes. In the line graph, the red line represents the results after multi-objective optimization, while the blue line denotes the monthly average water levels of Lake Superior after outlier removal.

\begin{figure}[H]
    \centering
    \includegraphics[width=0.55\linewidth]{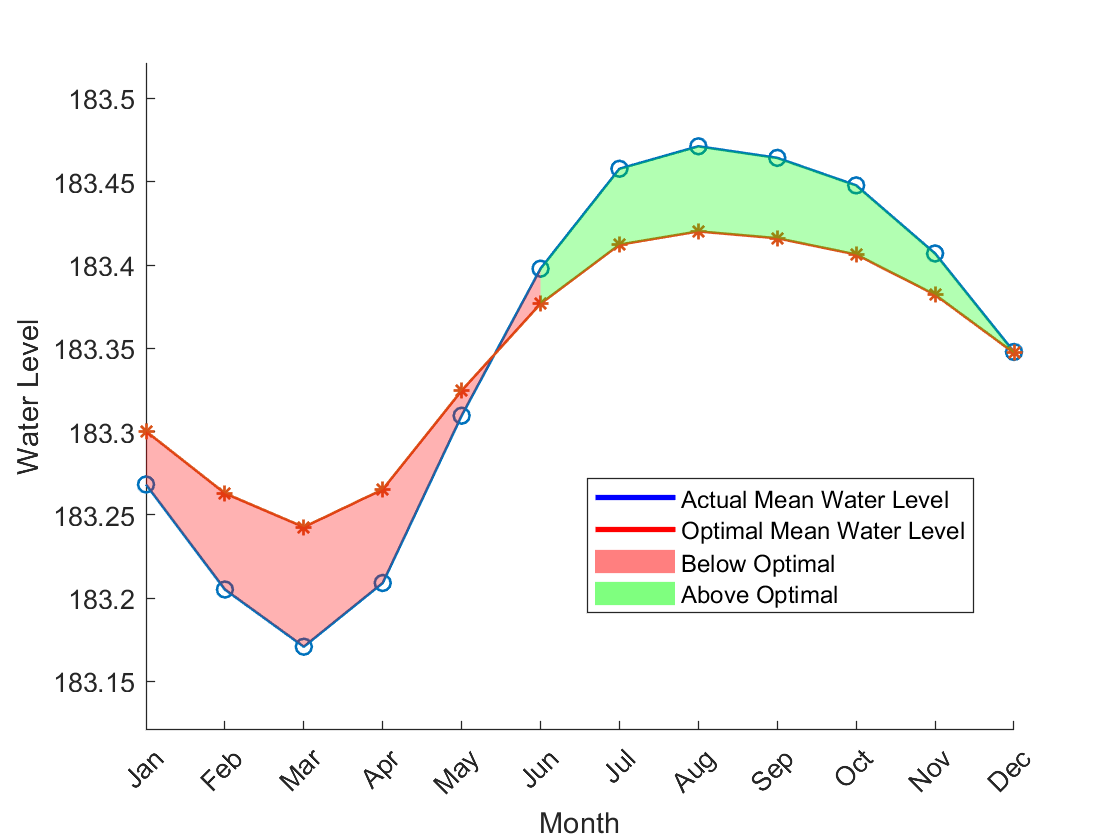}
    \caption{Optimal water level and actual mean level of Lake Superior}
    \label{fig:enter-label}
\end{figure}

The red fill indicates instances where actual levels fall below the optimal levels, whereas the green fill denotes the opposite, signifying instances where actual levels exceed the target levels. The color fill is applied on a monthly basis, reflecting our control metrics which are discretized into monthly values.

\section{Requirement 2 - DDE-MPC Approach to Predict Optimal Dam Regulation}
\subsection{Coefficient Determination}
Due to significant anomalies in the water levels and flow data from 2020 to 2023 when compared to the averages, which is likely to be influenced by unforeseen environmental factors during these years, the data for these four years is temporarily excluded from consideration. By inputting the upstream water level and flow data from 2009 to 2019 into the model and conducting data fitting, a linear function can be obtained. In this way, the coefficient $\beta$, $\gamma$, $\delta$, and $\Delta$ can be determined. The result of Linear fitting comes at:

\begin{figure}[H]
  \centering
      \subfigure[Erie-Niagara Relation]{\includegraphics[width=0.45\textwidth]{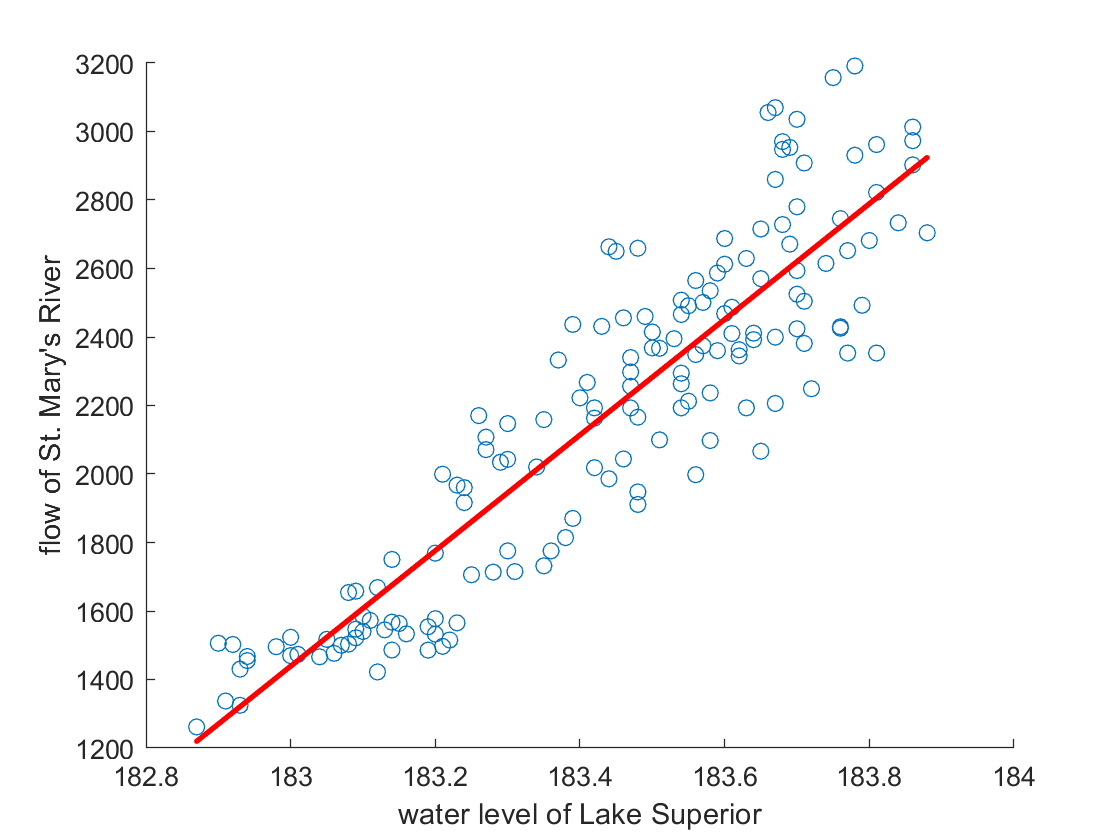}} 
    \subfigure[Clair-Detroit Relation]{\includegraphics[width=0.45\textwidth]{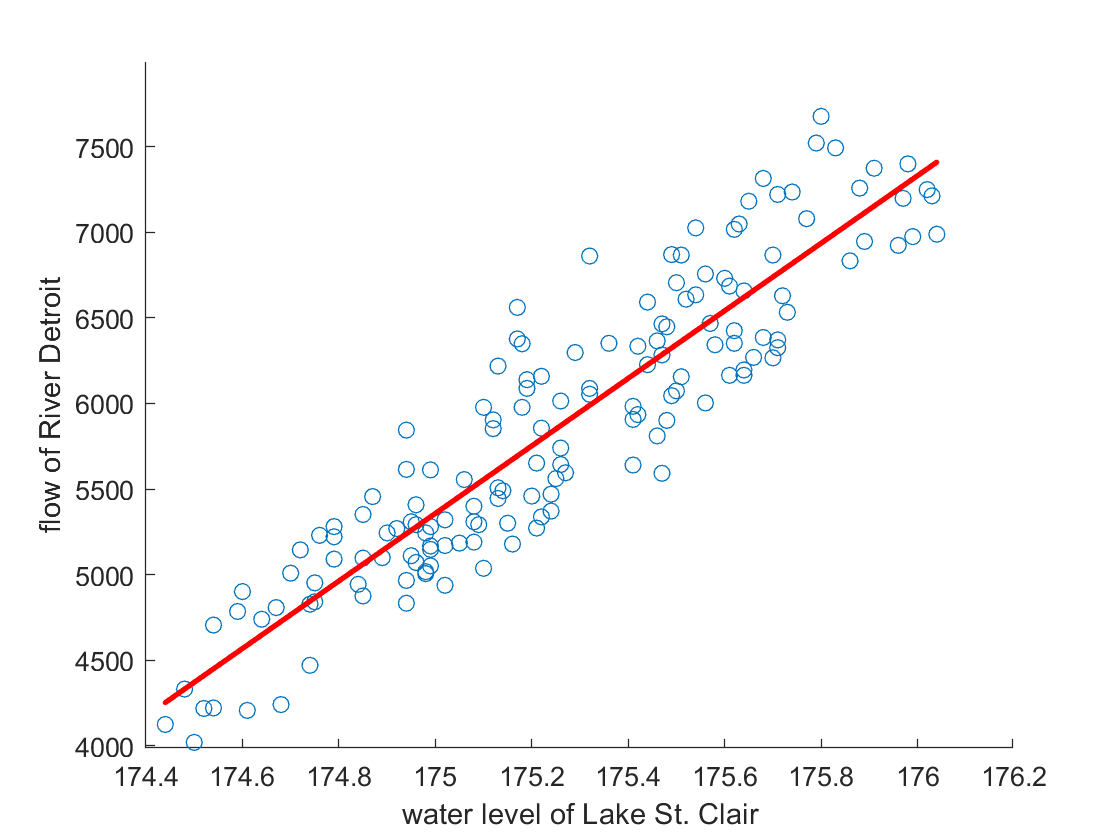}} 
    \subfigure[M\&H-Clair Relation]{\includegraphics[width=0.45\textwidth]{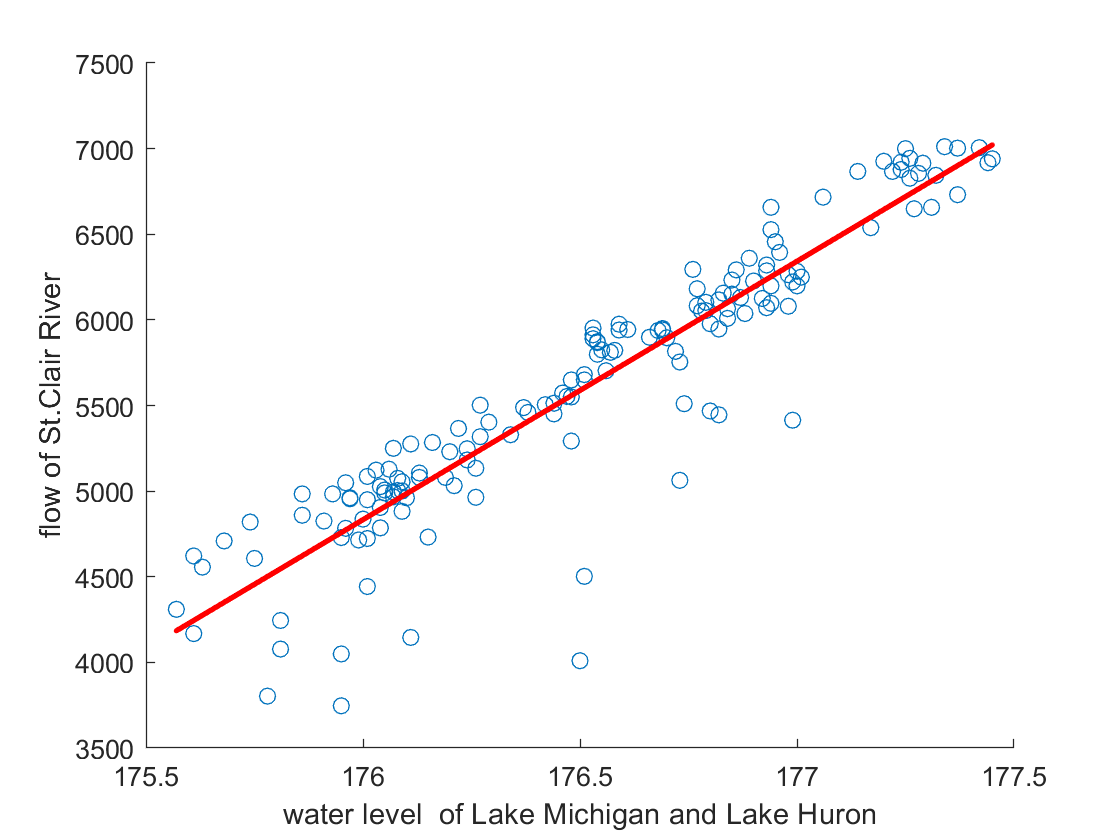}}
    \subfigure[Erie-Niagara Relation]{\includegraphics[width=0.45\textwidth]{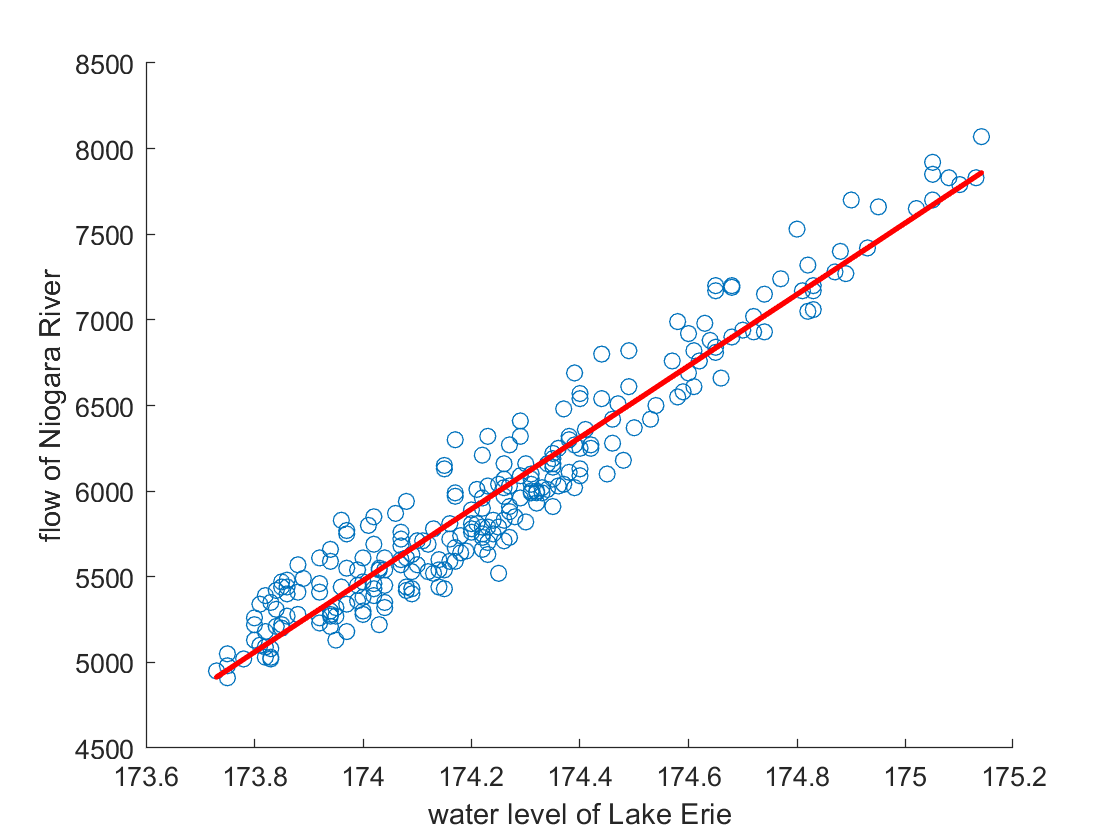}} 
  \caption{Linear Fitting of the relation between Lake water level and River flow}
  
    \label{fig:Linear fitting-Lakeriver}
\end{figure}

\begin{table}[htbp]
\begin{center}
\caption{Data fitting}
\begin{tabular}{ccc}
\toprule[2pt]
\multicolumn{1}{m{5.5cm}}{\centering Lake-River}
&\multicolumn{1}{m{4cm}}{\centering Slop$/10^3\cdot m^2s^{-1}$}&\multicolumn{1}{m{4cm}}{\centering Intercept$/10^5\cdot m^3s^{-1}$ }\\
\midrule
Superior-St.Mary & 1.69 & -3.0744 \\
  St.Clair-Detroit & 1.51 & -2.6076  \\
 Michigan and Huron-St.Clair & 1.97 & -3.3980 \\
Erie-Niagara & 2.09 & -3.5799 \\
\bottomrule[2pt]
\end{tabular}\label{tb:data fitting}
\end{center}
\end{table}

\subsection{DDE-MPC Based Prediction Model}
The model we apply in the requirement is an overall prediction system over Natural Indicator $\Delta A(t)$,and use the Indicator to signify the prospected lake water level in the coming months.
\subsubsection{Natural Indicators}
The overall water level of lake is generally in a circadian rhythm with yearly periodic. Apart from the inflow and outflow of the river that connects to it, Natural factors are also crucial to the nature-based-regulation of lake water level.
Major objective factors that contribute to the lake water level are climate, precipitation, evaporation, temperature, tide, seepage, altitude and moisture. With collection of the detailed data in the corresponding month, a general Natural indicator can be acquired.
The Natural indicator can also be indirectly calculated through the given data about the lake and rivers that connect to it. Therefore, the equation of lake B can be interpreted as:
$$
S\cdot\frac{\mathrm{d} B}{\mathrm{d} t} =In(t-1)-Out(t-1)+\Delta _{B}(t)
$$
The $In(t-1)$ and $Out(t-1)$ is computable through the Inflow and Outflow of river throughout the month. The equation is made in perspective of volume change. A table of Natural Indicators of Lake B is given as follows:
\begin{figure}[H]
  \centering
    \includegraphics[width=0.8\textwidth]{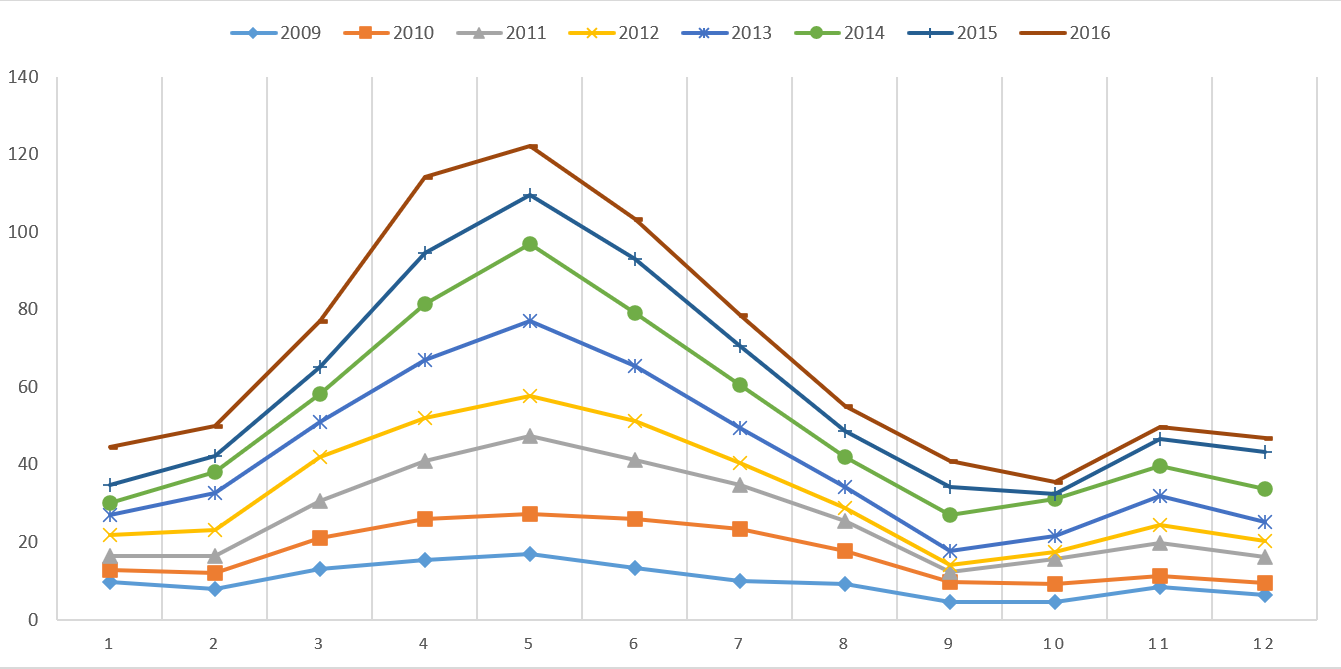}
  \caption{Lake B Natural Indicator over 8 years $\Delta_B$ (Unit:$10^8\cdot m^3$)}
    \label{fig:LakeBNatural}
\end{figure}

\subsubsection{DDE Model}
In the Model, a delayed-differential equation (DDE) perspective is applied to characterize the form of delayed impact river flow possesses on the corresponding lake. When generalizing the function, the current lake water level is estimated by the water flow inwards and outwards over the last month. Take Lake $B$ (Lake Huron and Michigan) as an example, the DDE can be written as:

$$
\frac{\mathrm{d} B}{\mathrm{d} t} =\beta _{11}a(t-1) -\beta _{12}b(t-1)+\Delta _{B}(t)\\
$$
The $\Delta$ here represents the Natural change of water. Therefore, a delayed response on the Great lake is added according to the alteration of dam, and the influence of the water regulation is postponed before problems are severe enough to be noticed. In reality, Plan 2014 aggravates the fluctuation of the Great Lakes, and the delayed response partially contributes to the adverse climate in 2017.

\subsection{Water Level Predictive Control Model}
We propose a control methodology predicated on Model Predictive Control (MPC), entitled the Water Level Predictive Control Model (WLPCM), which encompasses the following modules:

\textbf{Predictive Model}: At the heart of MPC lies the predictive model, necessitating a mathematical construct capable of forecasting the water level of a lake over a future period. This model integrates the effects of variables such as rainfall, evaporation, inflow, and outflow on water levels. For this purpose, we utilize an existing Delay Differential Equation (DDE)-based prediction model. Consequently, we employ monthly averages, devoid of outliers, as a predictive measure for the future, acknowledging that these estimates may diverge from actual outcomes.

\textbf{Optimization Solver}: The solution employs a Simulated Annealing algorithm to determine optimal dam flow rates over the forthcoming six months. Given the temporal lag of four months for upstream changes to manifest downstream, it is imperative to factor in the medium to long-term impacts, lest the model yields locally optimal solutions for the ensuing month at the expense of global optimality.

\textbf{Feedback Adjustment}: A pivotal component of MPC involves the collection of data on water levels and flow rates to calibrate the model and refine the control strategy. Consequently, even in the face of extreme weather conditions, our WLPCM is capable of providing a rational and reliable control response.

\textbf{Constraints}: Reasonable constraints are established based on the dam's maximum flow capacity and the safe range of water levels, ensuring the operational integrity and safety of the system.

\begin{algorithm}
\caption{Model Predictive Control Process}
\begin{algorithmic}[1]
\State \textbf{Input:} initial\_state, prediction\_horizon, control\_objectives
\State \textbf{Output:} control\_process\_results

\Procedure{MPC\_Process}{initial\_state, prediction\_horizon, control\_objectives}
    \State current\_state $\gets$ initial\_state
    \For{each control step to prediction\_horizon}
    
        \State \textit{// Define and solve the optimization problem}
        \State optimal\_control\_action $\gets$ Solve\_Optimization(Predict\_Model(current\_state, control\_action), control\_objectives, safety\_constraints)
        
        \State \textit{// Apply the optimal control action and update the current state}
        \State current\_state $\gets$ Apply\_Control(current\_state, optimal\_control\_action)
        
        \State \textit{// Collect the latest data for model correction and adjustment}
        \State Update\_Model\_with\_Latest\_Data(latest\_data)
        
        \State \textit{// Check if the control objectives are met or if there are any emergencies}
        \If{control objectives met or emergency occurs}
            \State \textbf{break}
        \EndIf
    \EndFor
    \State \Return control\_process\_results
\EndProcedure
\end{algorithmic}
\end{algorithm}

\begin{figure}[H]
  \centering
    \includegraphics[width=0.8\textwidth]{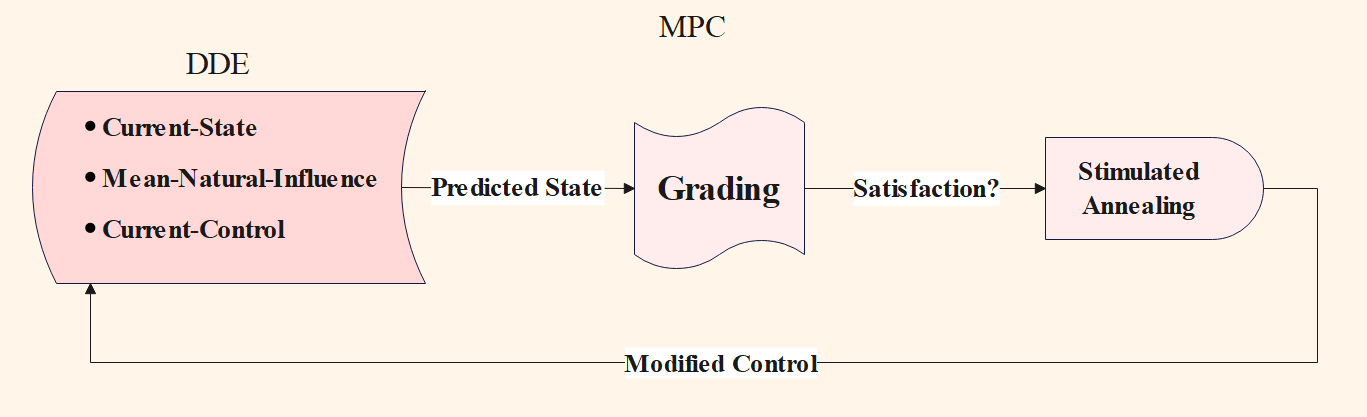}
  \caption{Model Predictive Control Process}
    \label{fig:mpc}
\end{figure}

\subsection{Prediction of future data}
We can apply the former Natural Indicator to postulate the newest Indicator $\hat{A_1}$ through linear regression. As DDE shows the postponed effect the dam regulation takes on the lake water level, we predict the following six months and use MPC to acquire an ideal dam regulation plan in the next month that reaches maximum of $G(a,e,R)$. Upon the completion of the first month, we are able to compute the Actual Natural Indicator $A_1$, which serves as a basis to forecast $\hat{A_2}$ for the ensuing month. The prediction lag retains 6 months, and advances in sync with the passage of time. As the year concludes, a comprehensive evaluation is made between the Actual Monthly Natural Indicator and the Projected Indicator, as well as the dam regulation $a[i]$ and $e[i]$ ($1\leq i \leq 12$), thus elucidating the algorithm's sensitivity. The main progress of the prediction is visualized in the subsequent picture.

\begin{figure}[H]
  \centering
    \includegraphics[width=0.9\textwidth]{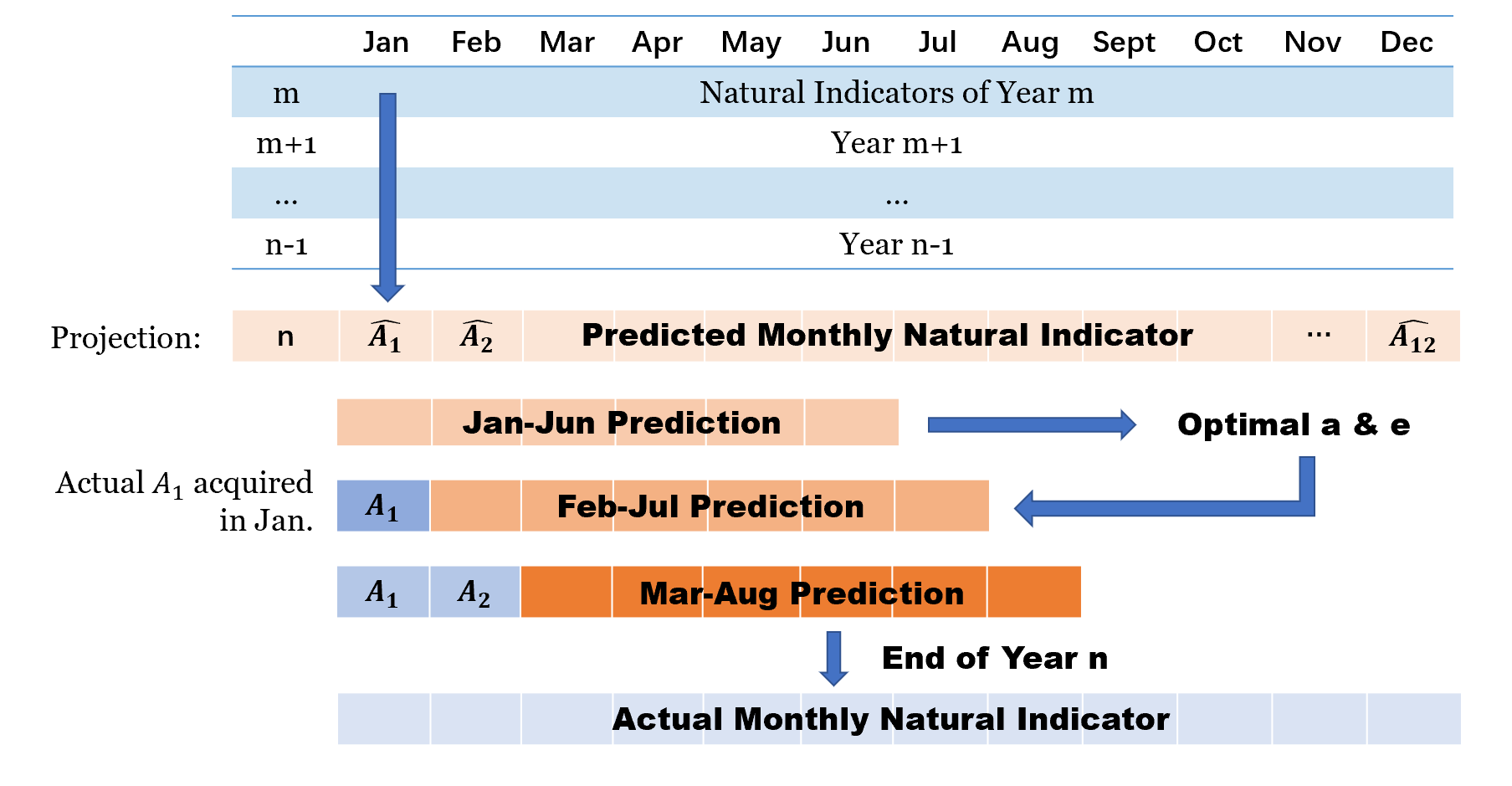}
  \caption{Model Predictive Control Process}
    \label{fig:ProDDEMPC}
\end{figure}

\section{Requirement 3 - Efficacy of WLPCM Model}
\subsection{Model Evaluation with Data in 2017}
\subsubsection{Evaluation Process}
To evaluate our WLPCM model, we employed another DDE Model distinct from the DDE integrated within the WLPCM. This particular DDE is equipped with knowledge pertaining to the natural weather conditions of the year 2017. Consequently, it can predict the actual values for the subsequent phase based on the current state and our control outputs. We utilized the simulation outcomes from this model to validate the performance of our model in the year 2017, thereby ensuring its accuracy and effectiveness.

\begin{figure}[H]
    \centering
    \includegraphics[width=0.95\linewidth]{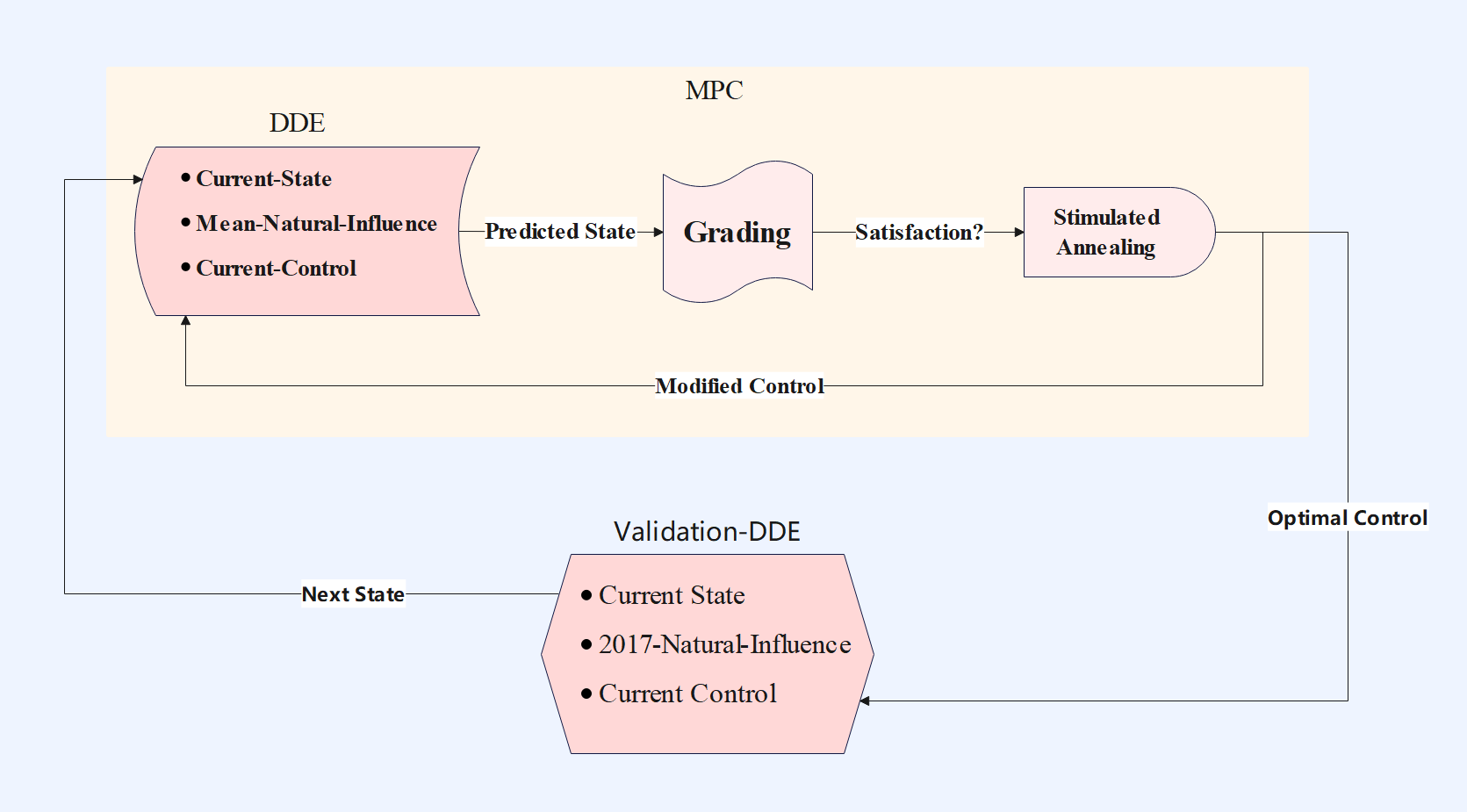}
    \caption{Evaluation Process}
\end{figure}

\subsubsection{Evaluation Results on Rivers}
Utilizing the WLPCM model, we executed control measures for the year 2017, resulting in the river water levels depicted in the figure below, which are compared with the actual flow.

\begin{figure}[H]
  \centering
    \subfigure[Actual River Water level in 2017]{\includegraphics[width=0.45\textwidth]{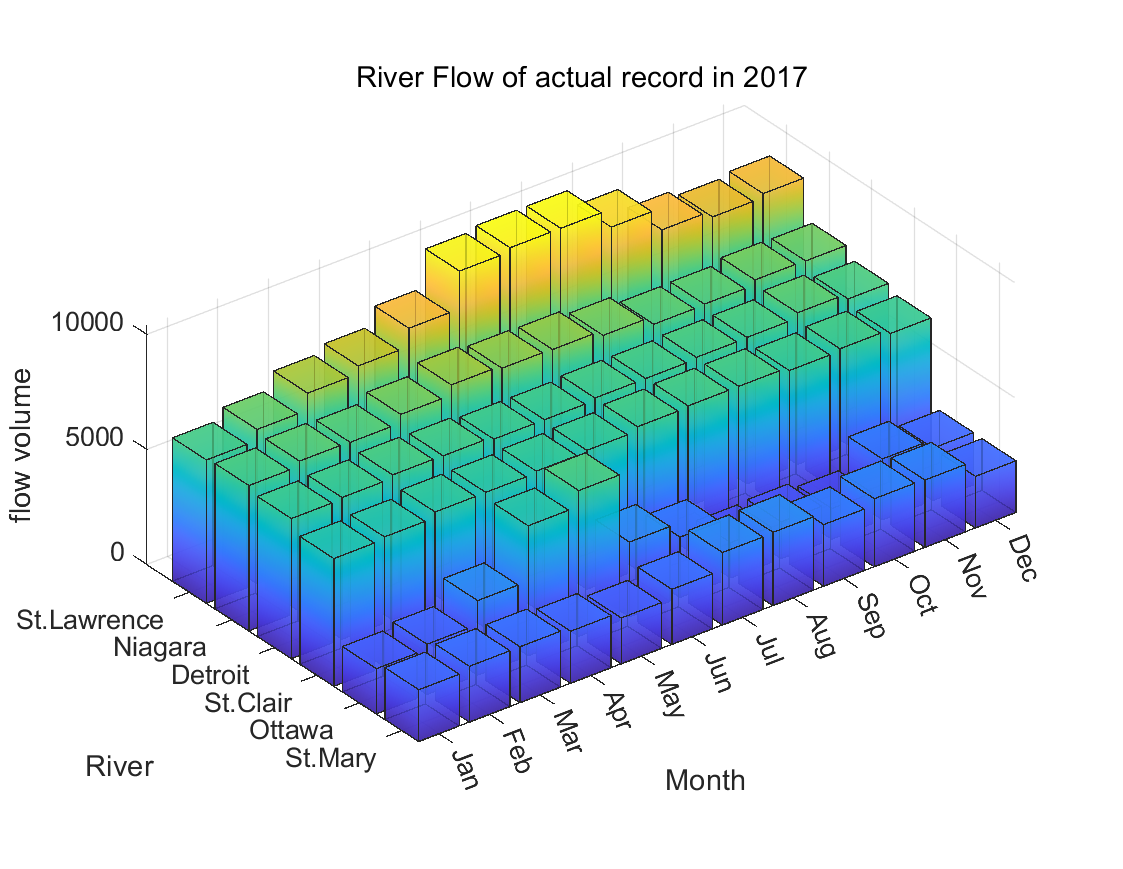}} 
    \subfigure[Optimized River Water level with our control]{\includegraphics[width=0.45\textwidth]{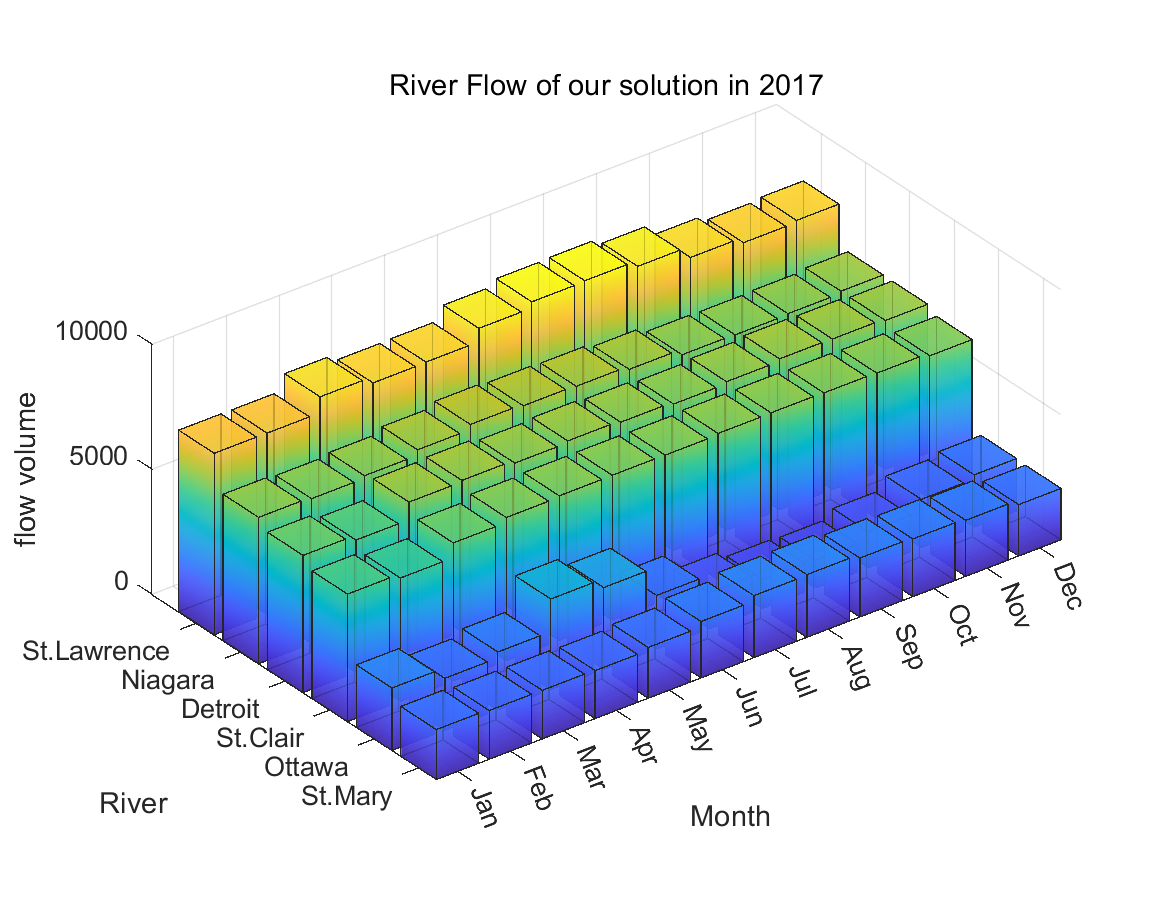}} 
  \caption{3D-bar Visualization of River Water Level}
\end{figure}

To visualize the actual scenarios of 2017 and the efficacy of our solutions, we defined a River Water Level Index. This index refers to the ratio of the river's monthly flow to the multi-year average, excluding outliers. This metric serves to quantify deviations from typical water levels, thereby facilitating a comparative analysis of hydrological conditions relative to historical norms.

$$
L_i(t) = \frac{flow_i(t)}{\overline{flow}_i(t)}
$$

\begin{figure}[H]
  \centering
    \subfigure[Actual River Water level in 2017]{\includegraphics[width=0.45\textwidth]{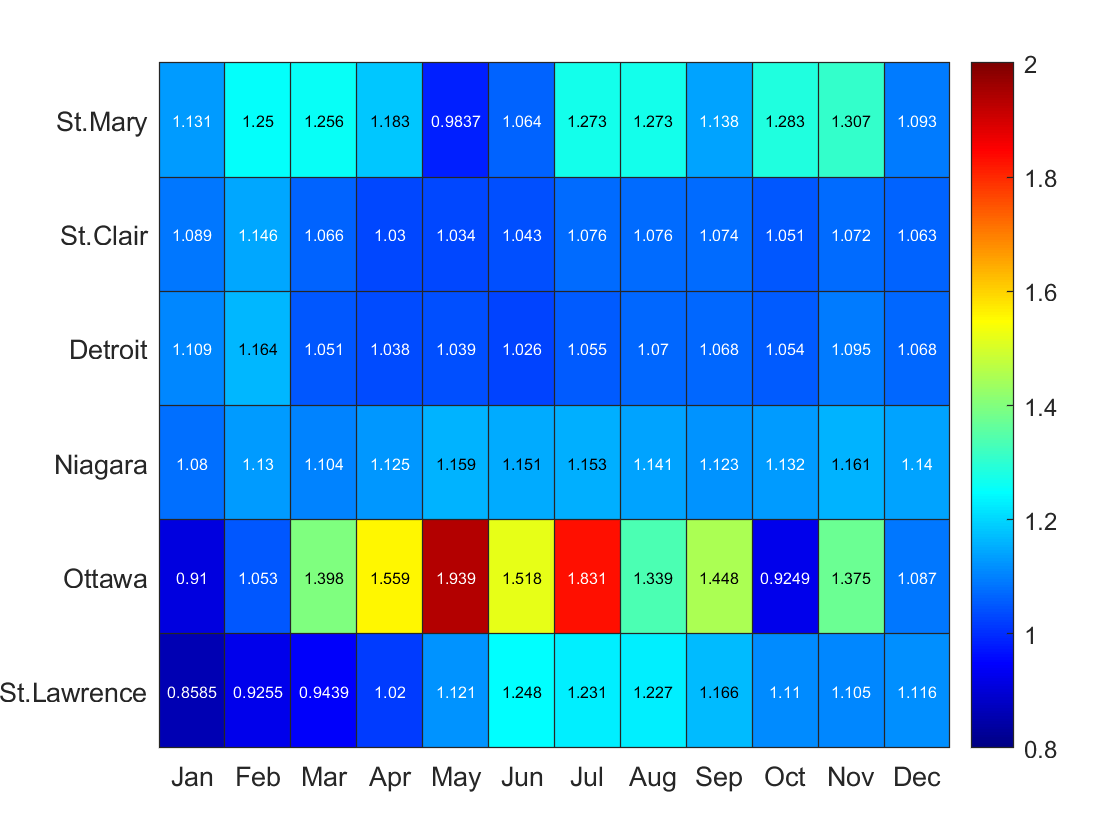}} 
    \subfigure[Optimized River Water level with our control]{\includegraphics[width=0.45\textwidth]{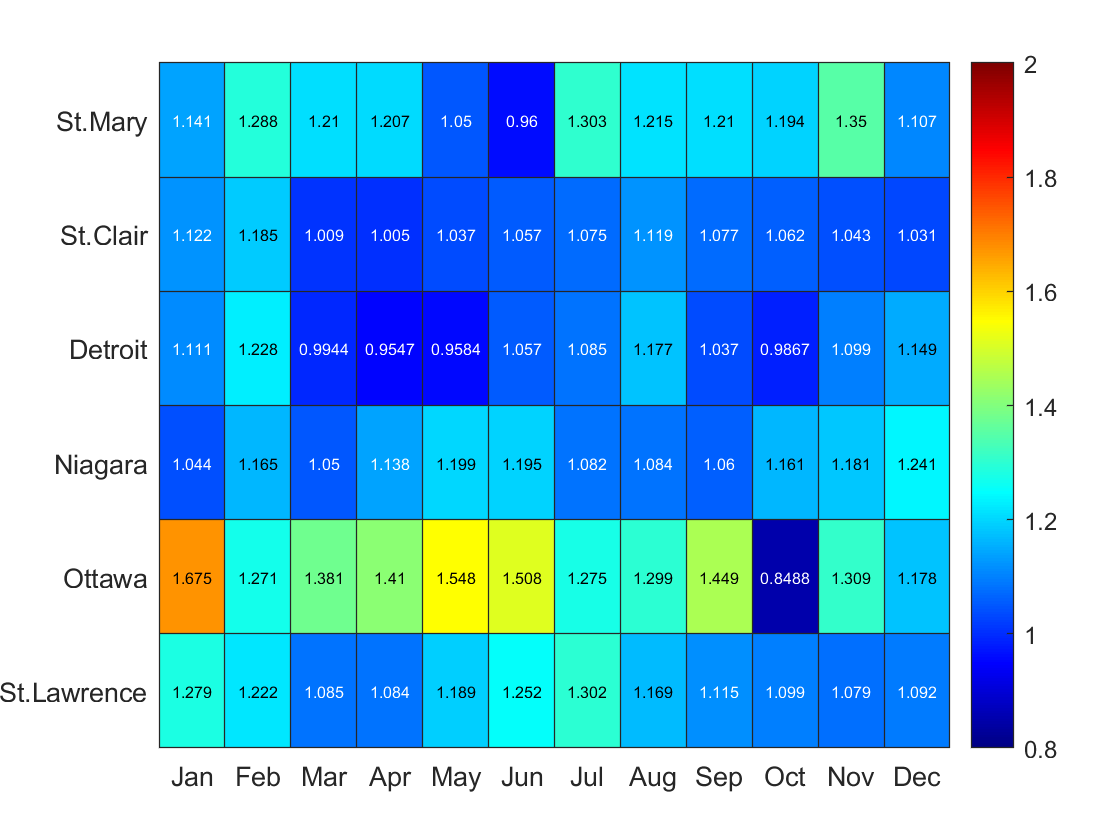}} 
  \caption{Heatmap of Rivers Water Level Index in 2017}
\end{figure}

\textbf{The anomalously high flow of the Ottawa River.}
In 2017, the Great Lakes region, along with the Ottawa River basin, experienced significant environmental challenges, most notably severe flooding. The flooding was a result of a combination of factors, including unusually high spring rainfall and the melting of a significant snowpack in the region. These conditions led to water levels rising to record heights in some areas, causing widespread damage to homes, infrastructure, and the environment.

\textbf{Under our control methodology, the substantial water flow is smoothly managed.}
In our heatmap, the water levels in the Ottawa are both high and stable, a condition that mitigates the risk of catastrophic flooding. This achievement is attributed to our six-month backward forecasting, which enables us to implement strategic long-term adjustments. By anticipating and counteracting potential short-term spikes in water levels, we effectively prevent disasters.

\subsubsection{Evaluation Results on Lakes}
For the water levels of lakes, given the relatively minor variations in lake water levels, we employ our defined Grading function for evaluation. It is noteworthy that within this context, $G_f$ represents the score for fluctuation magnitude, while $G_l$ denotes the score for water level.

\begin{figure}[H]
  \centering
    \subfigure[Actual 2017]{\includegraphics[width=0.45\textwidth]{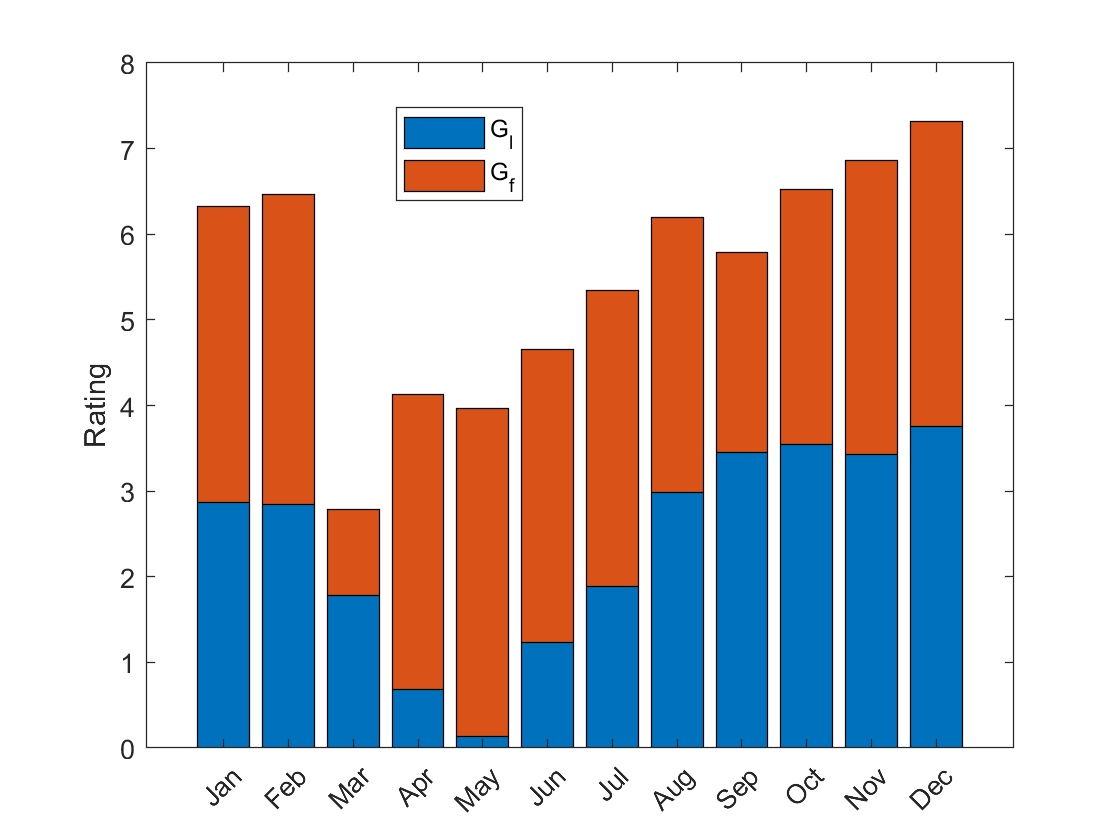}} 
    \subfigure[2017 with our control model]{\includegraphics[width=0.45\textwidth]{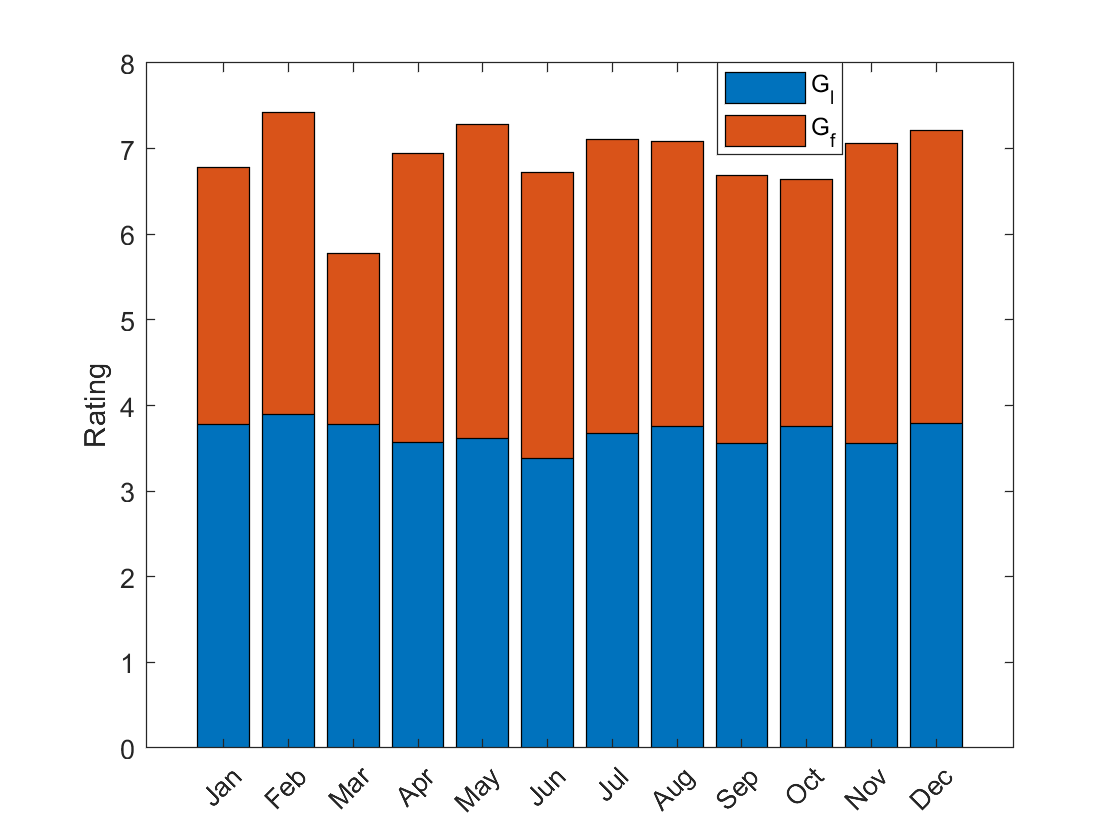}} 
    
  \caption{Ratings of Lakes Water Levels}
\end{figure}

For the annual score, we compare the yearly average, minimum value, and median of $G_f$ and $G_l$ respectively. This comparison allows us to evaluate the merits and demerits of our method against the actual results in terms of the annual lake water score.

\begin{table}[htbp]
\centering 
\caption{Comparison of Annual Water Level Scores ($G_l$)}
\begin{tabular}{ccc}
\toprule[2pt]
\textbf{Statistical Measure} & \textbf{Actual Records} & \textbf{WLPCM Method} \\
\midrule
Yearly Average & 2.38 & 3.68 \\
Minimum Value & 0.13 & 3.38 \\
Median & 2.85 & 3.71 \\
\bottomrule[2pt]
\end{tabular}
\end{table}

\begin{table}[htbp]
\centering 
\caption{Comparison of Annual Water Fluctuation  Scores ($G_f$)}
\begin{tabular}{ccc}
\toprule[2pt]
\textbf{Statistical Measure} & \textbf{Actual Records} & \textbf{WLPCM Method} \\
\midrule
Yearly Average & 3.14 & 3.22 \\
Minimum Value & 0.98 & 2.04 \\
Median & 3.44 & 3.36 \\
\bottomrule[2pt]
\end{tabular}
\end{table}

Based on the metrics listed above, whether it pertains to the degree of fluctuation or the levels of water, the effectiveness of our control surpasses that recorded in the year 2017 significantly.

\subsection{Sensitivity analysis}
\subsubsection{Methodology}
We take an Sobol Approach to analyze the sensitivity of the Lake-Water regulation. Through modifying the water flow of river $a$ and $e$, and analyzing the sensitivity by calculating the difference between the change value of the twelve-month water level variance and the change value of the water flow of the two rivers, an RMSE(The Root Mean Squared Error) outcome is acquired and can be used to indicate sensitivity of $a$ and $e$ on each lake.
$$
S_{i,\delta}=\frac{\sigma(G(a,b,R,i)-G(a-\delta,b,R,i)}{\delta}
$$

The RMSE estimation and Sobol approach boast following advantages:

\textbf{Error Magnitude:} RMSE gives a relatively straightforward interpretation of model error magnitude. A smaller RMSE indicates a better fit to the data, assuming the error metric is appropriate for the problem context.

\textbf{Scale Sensitivity:} RMSE values are in the same unit as the predicted and observed values, making them intuitively easier to understand. However, its sensitivity to the scale of the data means that RMSE values are more useful for comparing models on the same dataset rather than across datasets with different scales.

\textbf{Outlier Sensitivity:} Due to the squaring of each error term, RMSE is particularly sensitive to outliers. Large errors have a disproportionately large effect on RMSE, making it a useful measure when large errors are particularly undesirable.

\subsubsection{Results}
We modify $a$ and $e$ in the dam with the same Natural Indicator proportionally, generate the new lake water level and calculate RMSE index $S$, the result comes as follows:

\begin{figure}[H]
    \centering
    \includegraphics[width=0.6\linewidth]{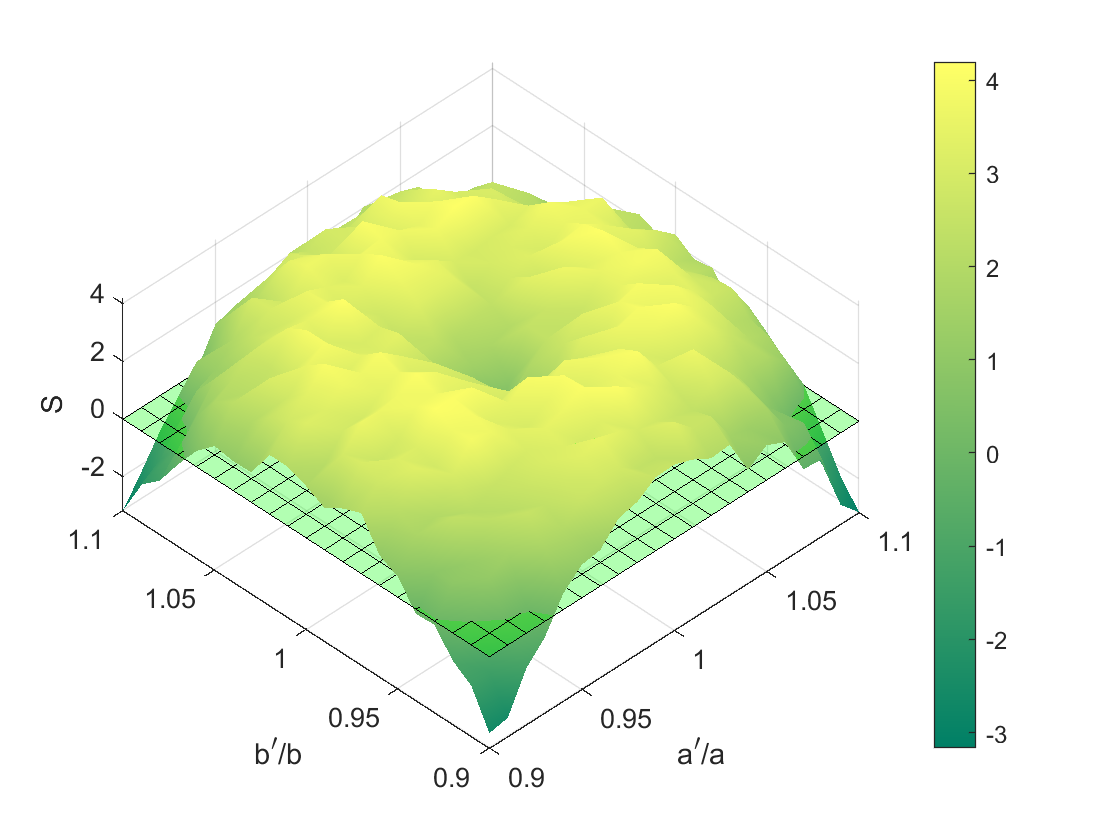}
    \caption{Sensitivity Test Result}
    \label{fig:p3-sensitivity}
\end{figure}
As illustrated in Figure~\ref{fig:p3-sensitivity}, our model demonstrates reasonable adjustments in its outputs in response to changes in input conditions. Hence, our model exhibits exceptional performance in this sensitivity test.
\section{Requirement 4 - Sensitivity Analysis}
\subsection{Modelling}

With given situation $M$, the sensitivity can be defined through fitness function $G$ as:
$$
S_{i,\delta}=\frac{\sigma(G(M+\delta,R,i)-G(M-\delta,R,i))}{\delta}
$$
where $\delta$ signals alteration of certain change of natural circumstances. The total Sensitivity of the system is defined as:
$$
S = S_{Rain}+S_{Ice}+S_{Snow}
$$
Subsequent discourse is dedicated to elucidating the sensitivity analysis through the examination of three distinct environmental alterations, encapsulated as follows:
\begin{enumerate}
\item \textbf{Precipitation}
introduces a modification in the Natural Indicators $\Delta$, functioning as an additional input to the system cycle. The alteration $\delta$ pertains to the modulation of Natural Indicators, with the sensitivity analysis aimed at elucidating the impact of these indicators on the system's dynamics.

\item \textbf{Ice Clog} 
results in a diminished water flow during the winter months, subsequently experiencing an elevation in the spring. Within the confines of our model, $\delta$ gives a decrement in river flow by $\delta$ in January and an increment by $2\delta$ in March. The interrelation between the degree of ice clog $\delta$ and water flow is graphically represented as follows:

\begin{figure}[H]
    \centering
    \includegraphics[width=0.6\linewidth]{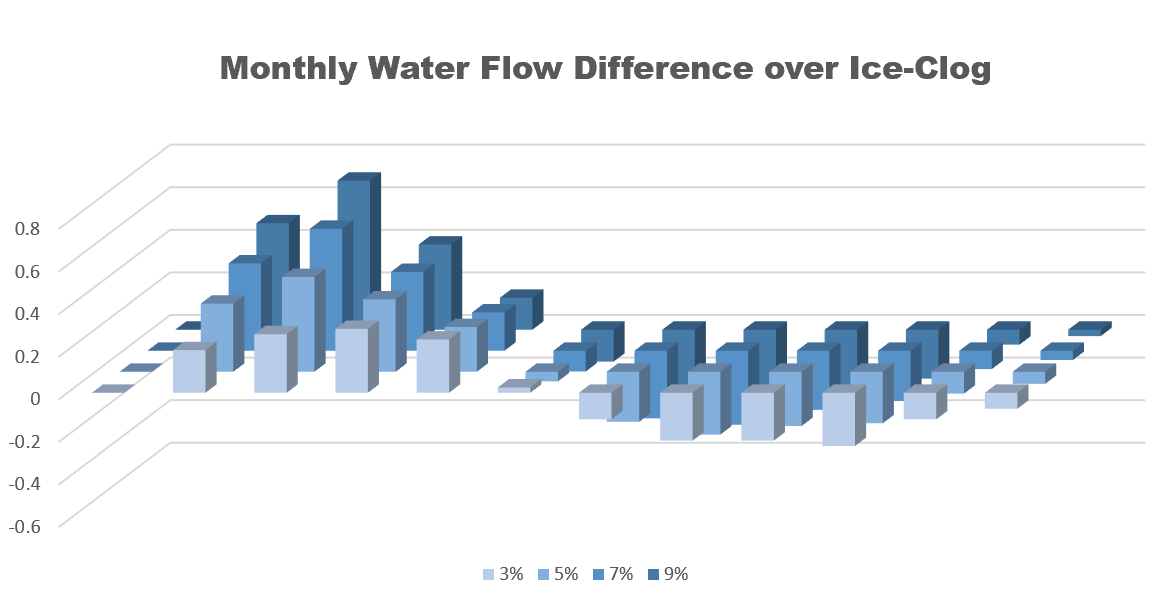}
    \caption{Relation between ice clog degree $\delta$ and water flow}
    \label{fig:Ice-clog}
\end{figure}

\item \textbf{Snow pack} during the winter months increased meltwater in the spring, thereby elevating the water levels in Superior Lake and exerting a cascading effect throughout the Great Lakes system. $\delta$ here indicates the March water level rise in Superior Lake.
\end{enumerate}
\subsection{Outcome}
We analyze the influence $\delta$ has on the scoring of level and fluctuation of $G$ and reach the following outcome: 

\textbf{Precipitation} ($\Delta_{Rain}$): A \textbf{3\%} increase in precipitation results in a \textbf{5.153\%} rise in the system's water level, indicating high sensitivity to rainfall changes.

\textbf{Ice Clogs} ($\Delta_{Ice}$): The introduction of a \textbf{3\%} decrease in water flow due to ice clogs in winter, followed by a \textbf{3\%} increase in spring, results in a \textbf{4.426\%} overall fluctuation in water levels. This showcases a moderate sensitivity to ice clog conditions.

\textbf{Snow Pack} ($\Delta_{Snow}$): A \textbf{3\%} increase in snow pack leads to a \textbf{2.707\%} increase in water levels during the meltwater phase in spring, signifying a very high sensitivity to snow accumulation.

The system exhibits the highest sensitivity to ice clog changes, followed by precipitation and snow pack. This indicates that ice clog dynamics play a crucial role in the Great Lakes system.
The differential sensitivity to these factors underscores the importance of considering seasonal variations and climate patterns in managing water resources.

The outcome of this sensitivity analysis suggests prioritizing the management of ice clog and precipitation to mitigate potential risks associated with water level fluctuations.
Strategies may include enhancing water storage capacity in anticipation of meltwater influx and implementing measures to manage ice formation and melting to stabilize winter and spring water flows.

\section{Requirement 5 - Ontario-based Pragmatic Problem}
\subsection{Needs or wants of key stakeholders}
To narrow our extensive analysis exclusively to the stakeholders and factors influencing Lake Ontario, we meticulously examined the demands of the stakeholders associated with Lake Ontario.

For each category of stakeholders, an analysis was conducted to identify their respective needs or desires, the rationale they employ to justify their positions, the foundational values and beliefs underpinning their stances, and any anticipations they harbor concerning the resolution process. A synthesized overview of the principal stakeholders' needs and desires is presented in the table below.\cite{furber2016conflict}

\begin{table}[H]
\centering 
\caption{Needs or wants of key stakeholders}
\begin{tabular}{p{5cm}p{10cm}}
\toprule[2pt]
\textbf{Stakeholder} & \textbf{Needs/Wants} \\
\midrule
Shipping and Navigation & Minimum water levels \\
Hydropower & Consistent flow \\
Recreational Fishing, Boating and Tourism & Abundant fish, extended boating season through higher water levels at end of season \\
Opposed South Shore & Consistent water levels \\
Environmental Groups & seasonal high high-water levels and low low-water levels on Lake Ontario \\
\bottomrule[2pt]
\end{tabular}
\end{table}

\subsection{Lake-river relation}
Given the multitude of water level demands for Lake Ontario, we will specifically focus on the water level management situation for Lake Ontario. The change of the water level of Lake Ontario is still determined by natural factors and the past months water flow data inwards and outwards, which can be interpreted as:
$$
\frac{\mathrm{d} E}{\mathrm{d} t} =\eta _{11}d(t-1) -\eta _{12}e(t-1)+\Delta _{E}(t)
$$
Moreover, as the livelihoods of the Montreal residents are directly influenced by the water dynamics accumulated in the city of Montreal, and the city water flow is contributed by the inward volume of water flowing through the Moses-Saunders Dam and the Ottawa river, and the outward water to the St. Lawrence River, the water flow in Montreal can be interpreted as:
$$
F_{St.Lawrence}-F_{Ottawa}-F_{e} = \Delta_{Montreal} = \Delta_{Nature} + \Delta_{Residents} + \Delta_{Retention} 
$$

\subsection{Fitness Function}

\subsubsection{Restrain on Lake}
To mitigate the impact of potential flooding from Lake Ontario on the residents of Montreal, we need to assess the risk of flooding in Lake Ontario. Therefore, apart from the $G_{L}$ and $G_{F}$, the grading of the flooding potential of Lake Ontario can be defined as:
$$
G_{D}=\left\{
\begin{aligned}
&0; & H<H^\#\\
&-4\sqrt{\frac{H-H^\#}{F_\sigma}}; & H^\#<H<H_{highest} \\
&-4; & H>H_{highest}
\end{aligned}
\right.
$$
Here, $H^\#$ represents the warning flood level, which is determined by the average water level of Lake Ontario and the std.deviation. $H^\#= F_{Ave}+F_{\sigma}$. And $H_{highest}$ is the highest water level of Lake Ontario, defined by
$H_{highest}=F_{Average}+2F_{std.deviation}$.

\begin{figure}[H]
  \centering
    \includegraphics[width=0.5\textwidth]{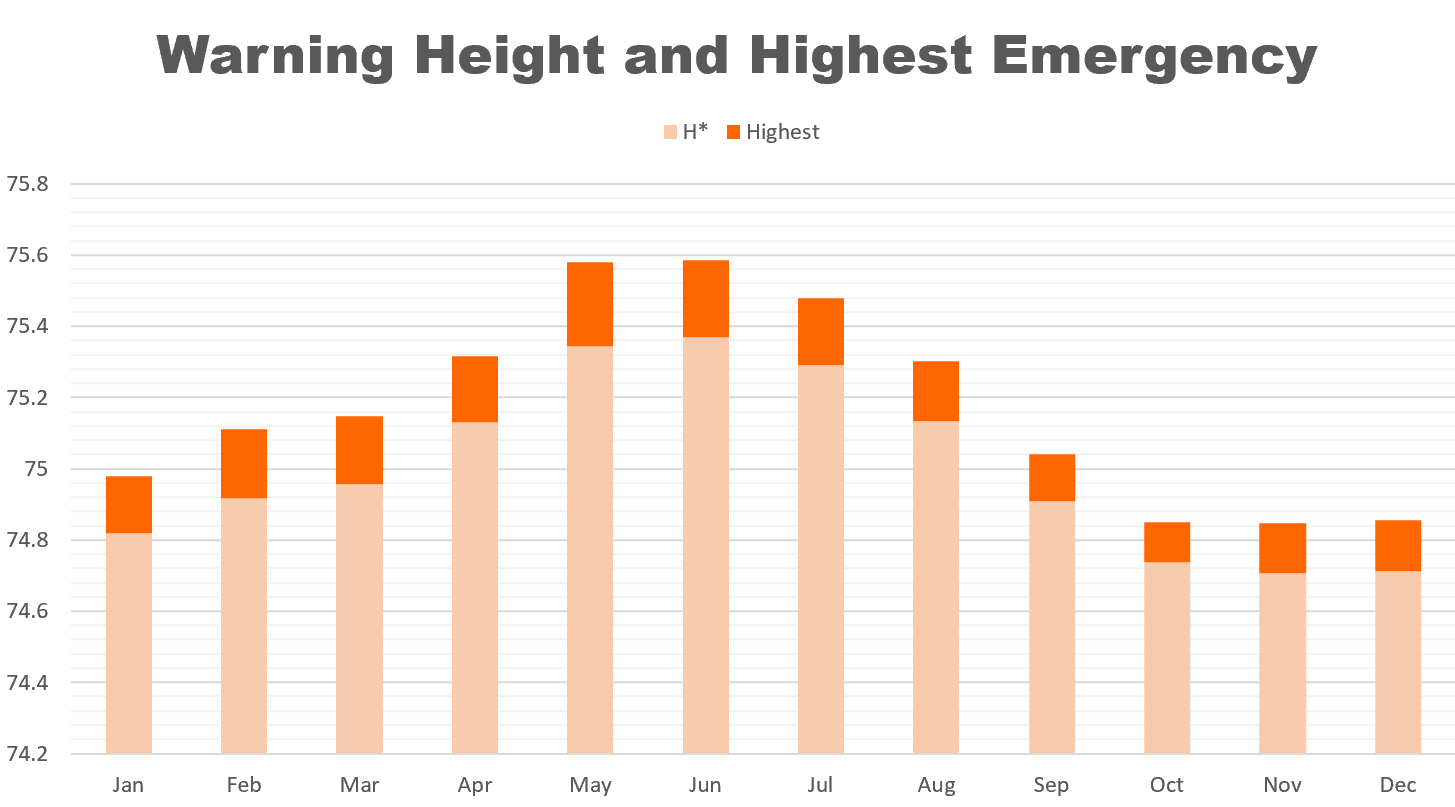}
  \caption{The warning height and the highest emergency of Lake Ontario}
    \label{fig:warning}
\end{figure}

\subsubsection{Restrain on River}
\textbf{Grading on River:} Similar to the lake surface water levels, we can define the grades of river both in water surface flow and river fluctuation  degree.

The grading of water flow can be defined as
$$
G_{L'}=\left\{
\begin{aligned}
&2+200(F-F^*);   & High-level-required \\
&4-400|F-F^*|; &  Medium-level-required \\
&2-200(F-F^*); &  Low-level-required
\end{aligned}
\right.
$$
and $F^*$ indicates the average water flow in this situation.

The grading of river fluctuation can be defined as
$$
G_{F'}=\left\{
\begin{aligned}
&2+1/80(\sigma - \hat{\sigma})  ; & High-fluctuation-required \\
&4-1/160|\sigma - \hat{\sigma}|  ; & Medium-fluctuation-required \\
&2-1/80(\sigma - \hat{\sigma})  ; & Low-fluctuation-required
\end{aligned}
\right.
$$
Here $\sigma$ and $\hat\sigma$ stands the newest standard deviation and original standard deviation:
$$
\sigma = \sqrt{\frac{\sum\limits_{i=1}^{12}(F-F^*)^2}{12}}
;\quad
\hat\sigma = \sqrt{\frac{\sum\limits_{i=1}^{12}(\hat{F}-\hat{F^*})^2}{12}}
$$

\textbf{Water retention in Montreal:} To prevent the increasing water levels from having potential adverse impacts on residents in Montreal, the hope for no water retention in Montreal arises. This involves safeguarding infrastructure, ensuring resident safety and well-being, and maintaining hydro power to foster a secure and resilient community environment. The goal is to uphold optimal water level, hoping that the water retention in Montreal reaches minimum. The water retention in Montreal can be defined as 
$$
G_{M}=-(F_{St.Lawrence}-F_{Ottawa}-F_{e}-\Delta_{Nature} - \Delta_{Residents})^2
$$

\subsection{Methodology}
In our approach to solving for the water levels of Lake Ontario, we continue to employ Model Predictive Control (MPC). Unlike before, our current Delay Differential Equation (DDE) model offers more accurate predictions of the water levels in Lake Ontario as well as the flow rates in the Ottawa and St. Lawrence Rivers. The grading module, at present, only examines the water levels of Lake Ontario and the flow rates of the St. Lawrence River; hence, we have exclusively considered the stakeholders of Lake Ontario in our analysis. We utilize a revised Water Level Prediction and Control Model (WLPCM) for forecasting, along with an updated Validation-DDE and Grading for assessment, ensuring a comprehensive and refined evaluation.

\subsection{Result on 2017}
\subsubsection{Records of 2017}
In 2017, the Lake Ontario-St. Lawrence River basin experienced unprecedented spring precipitation and snowpack runoff, leading to historic high water levels and overwhelming regulatory flood mitigation efforts. The governing 2014 Plan by the International Joint Commission aimed to balance water levels for environmental and human use but faced challenges under these extreme conditions. The resulting severe flooding triggered a reevaluation of this plan, highlighting the difficulties of managing such a vast, cross-border water system amid climate change. It stressed the importance of adaptive management, international cooperation, and enhanced research for future resilience.

\subsubsection{Results of Our Control}
\begin{figure}[H]
  \centering
    \subfigure[Water Level of Lake Ontario]{\includegraphics[width=0.45\textwidth]{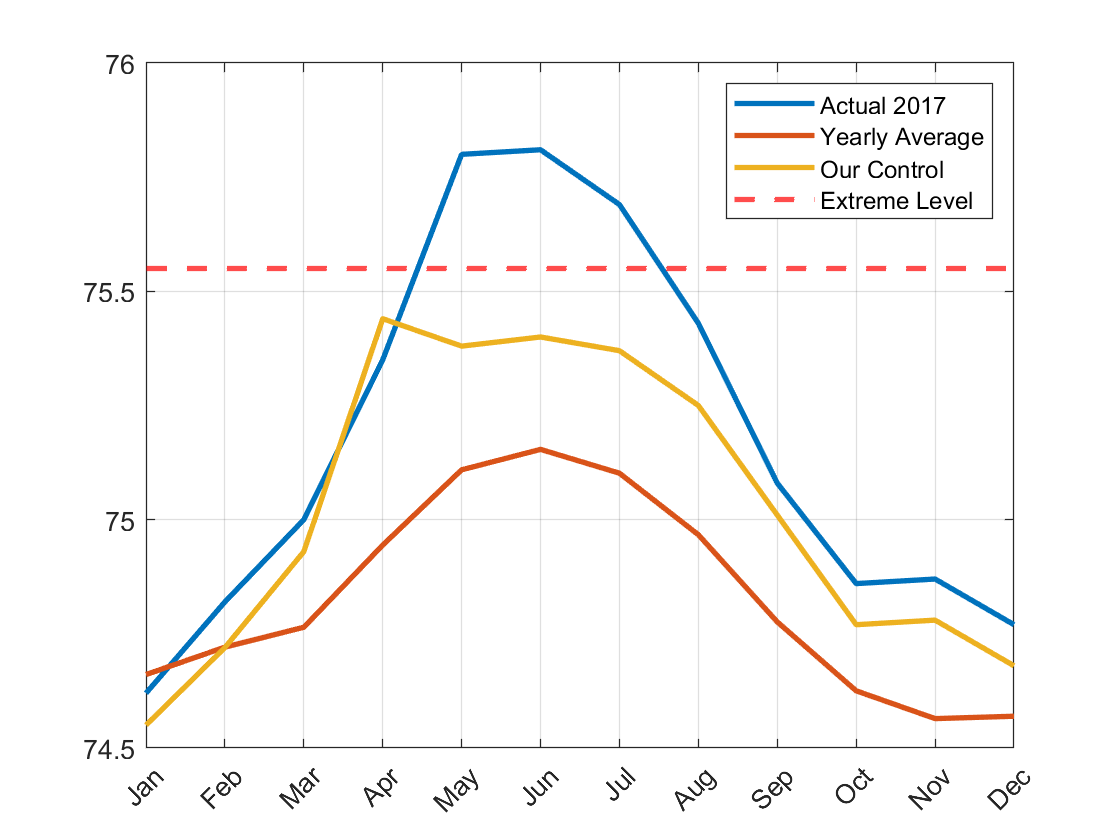}} 
    \subfigure[Flow of River St.Lawrence]{\includegraphics[width=0.45\textwidth]{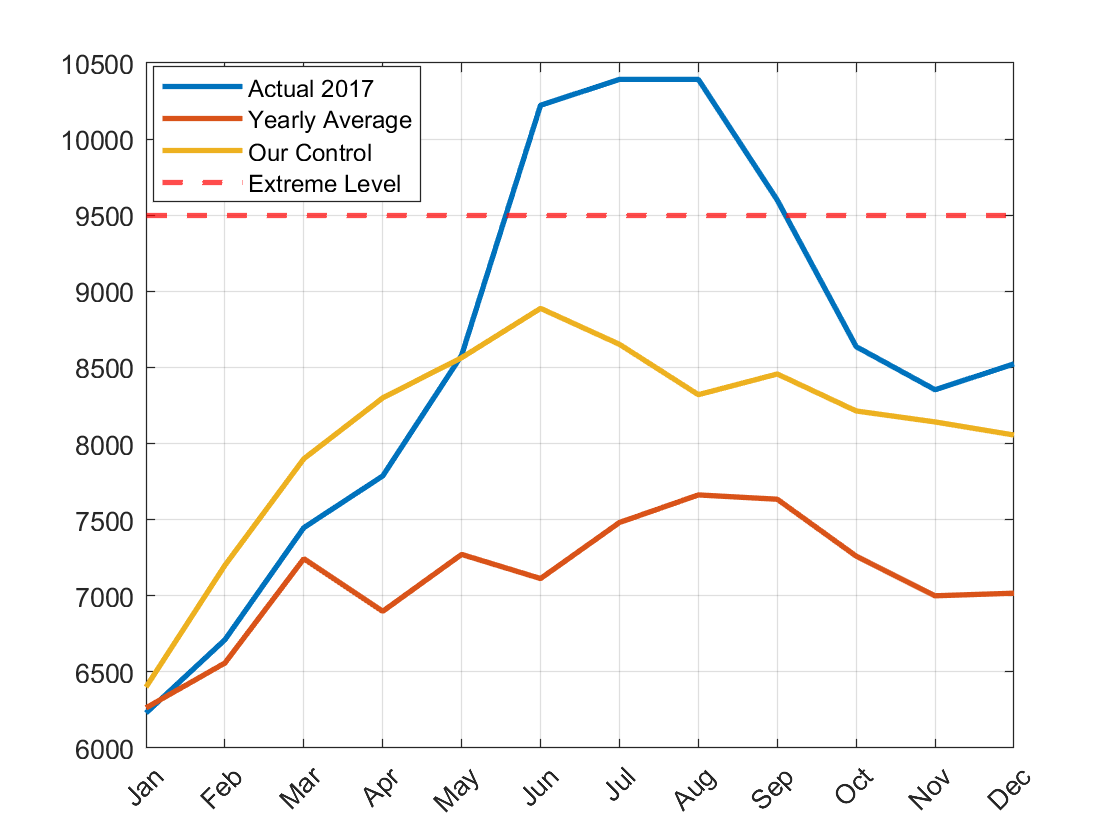}}
    
  \caption{Flow of River St.Lawrence}
    \label{fig:prob5-control compare}
\end{figure}

Based on the model we previously established, we conducted simulations using data from the year 2017. Figure~\ref{fig:prob5-control compare} illustrates the water levels of Lake Ontario and the flow rates of the St. Lawrence River as generated by our control model.

Referencing reports on flooding from that year, we delineated the extreme water levels as depicted in the figure provided. Compared to the Plan 2014 utilized that year, our control method significantly reduced peak water levels, thereby preventing the occurrence of floods. Moreover, despite these interventions, the fluctuation of water levels in Lake Ontario remained above the historical average, ensuring no compromise on the conservation of biodiversity. Consequently, in terms of managing the water levels of Lake Ontario, our control method demonstrates comprehensive superiority over Plan 2014.

\section{Model evaluation and further discussion}
\subsection{Strength}
\begin{itemize}
\item \textbf{Optimal over the long term}\\
Our model incorporates a Delay Differential Equation model to forecast future natural conditions, optimizing for the best solution over a six-month horizon. This approach enables the implementation of long-term optimal control measures, thereby mitigating the occurrence of extreme disasters.
\item \textbf{Capable of wide application}\\
By inputting historical data of the river-lake system to be controlled, along with the needs of stakeholders, our model can provide long-term optimal solutions. Thus, it is applicable to any location within the Great Lakes or other river-lake systems around the world.

\item \textbf{Robust}\\
Our WLPCM model incorporates a feedback mechanism that accepts data generated from the natural world after executing control outputs, refining future forecasts for more accurate control. With the intensification of climate change and the expected increase in extreme weather events, our model's ability to maintain precise control under unpredictable conditions is indispensable. This robustness is essential in responding to the challenges posed by climate change.
\end{itemize}

\subsection{Weakness}
\begin{itemize}
\item \textbf{Demands significant computational capacity}\\
Each control iteration in our model involves forecasting and optimization for the future, necessitating substantial computational effort. However, with the advancement of modern computational hardware, such as GPUs, our model is capable of rapidly responding to external data inputs and delivering optimal control solutions.
\end{itemize}


\nocite{clamen2018plan}
\nocite{koops2015lake}
\nocite{Minallah2023ControlsOV}
\nocite{white2008natural}
\nocite{international2018observed}
\nocite{gronewold2013coasts}
\nocite{razavi2014evaluation}
\bibliography{reference}




\newpage
\textbf{To}: {The International Joint Commission (IJC)}

\textbf{From}:{Team \#2401540}

\textbf{Date}:{Feb.5, 2024}

\textbf{Subject}: {WLPCM model: A DDE-MPC approach to optimize management of the Great Lakes}

Our team is writing to present a new predictive model to manage Great Lakes based on the MPC algorithm, we name it as the WLPCM model based on DDE-MPC approach.

The essence of WLPCM lies in the combination of a \textbf{predictive model} and \textbf{simulated annealing} through grading of the fluctuation and the average water level throughout the year according to stakeholders' desires, then generating the optimal control of Compensating Works and Moses-Saunders Dam. We predict the impact of the dam control over the next six months, then employ the best recipe in our prediction, and continue to upgrade the postulation with the passage of time.

The \textbf{Delayed Differential Equations} explain the relationship between lake level and river flow over the previous month and simulate the linear relations to simplify calculation. \textbf{Model Predictive Control} serves as the methodology to evaluate consistency between dam control and stakeholder desire via fitness function that conducts grading process. This methodology not only enhances the accuracy of our environmental forecasts but also significantly improves our ability to respond dynamically to emerging challenges, ensuring sustainable water quality and ecosystem health.

The major target of our prediction is the \textbf{Natural Indicator} representing the influence of nature on the Great Lakes system. We applied the data based on climate (i.e. evaporation and seepage), river flow and water level over 1964-2016 before Plan 2014 has been employed to calculate and validate the historical Natural Indicator and apply DDE-MPC to predict future Indicator to carry out the relevant calculation.  The simulation results of our model based on 2017 data indicate that, in terms of water level satisfaction, our model achieved a \textbf{54.62\%} improvement over Plan 2014, and in fluctuation satisfaction, it improved by \textbf{2.55\%}.

In view of the specific approach focused on the Lake Ontario, St. Lawrence River and Montreal city, the fitness function is an overall balancing determined by factors like water level, fluctuation and flood prevention, and could greatly avoid the overflowing of Lake Ontario compared with the current Plan 2014 that renders extremely high-water level and the flood related to it because the plan stick to a 'higher-high, lower-low' approach to maintain diversity of environment.

In conclusion, the WLPCM model represents a strategic and proactive approach to addressing the environmental, social, and economic challenges facing our shared waters. We believe that this model not only aligns with but actively supports the IJC's mission and objectives. We are eager to engage with the Commission to explore how this program can complement and enhance ongoing efforts under the auspices of the IJC.

We look forward to the opportunity to discuss this program further and to explore potential avenues for collaboration.

Sincerely,

Team \#2401540






\AImatter

\begin{ReportAiUse}{9}
\bibitem{AI1}
OpenAI ChatGPT 4\\
Query1: Introduce the basic format of a memo
\\
Output: A memo, short for memorandum, is a type of document used in business and official communications to convey information within an organization. It is designed to be a quick and efficient way to communicate policies, procedures, announcements, or any other internal matters. The basic format of a memo typically includes the following elements:

Header: At the top of the document, the header clearly indicates that the document is a memo. It typically includes fields for the following:

To: The recipient(s) of the memo.

From: The sender of the memo.

Date: The date when the memo was written.

Subject: A brief summary of the memo’s purpose.

Opening: The opening segment briefly states the purpose of the memo. It should be concise and to the point, providing a clear reason why the memo is being sent.

Context or Background: This section provides any necessary background information that the recipient(s) need to understand the memo's purpose and content. It sets the stage for the discussion or the announcement that follows.

Discussion: The main body of the memo elaborates on the information or the announcement. It details the subject mentioned in the opening, providing all necessary details, explanations, and justifications. The discussion should be clear and structured logically.

Action Items or Conclusion: This part outlines any actions that need to be taken, deadlines, or it summarizes the memo’s key points. If specific responses or actions are required from the recipients, those should be clearly stated, including any deadlines.

Closing: While not always necessary, a closing statement can be included for courtesy or to offer additional assistance or contact information for follow-up questions.

Attachments (if applicable): If the memo references other documents or includes additional information as attachments, this should be noted at the end of the memo with a brief mention of each attachment.

The tone and style of a memo can vary depending on the audience and the message. However, memos are generally written in a concise, formal tone to ensure clarity and professionalism.

\bibitem{AI2}
OpenAI ChatGPT Classic\\
Query1: introduce the Conditions and Regulated Outflows of Lake Ontario-
St. Lawrence River in 2017\\
Output: In 2017, the management of Lake Ontario-St. Lawrence River water levels and flows underwent significant scrutiny and adaptation, particularly in response to the record-high water levels experienced in the Lake Ontario and St. Lawrence River system. The conditions and regulated outflows during this period were governed by a combination of climatic factors and regulatory frameworks designed to balance the diverse interests affected by the water levels.

Climatic Conditions
2017 was marked by exceptionally high precipitation in the Lake Ontario-St. Lawrence River basin during the spring months, contributing to record-breaking water levels in Lake Ontario. This was compounded by the runoff from a substantial snowpack in the upper Great Lakes and the Ottawa River basin, which flows into the St. Lawrence River upstream of Montreal. These natural conditions exceeded the capacity of any regulatory plan to mitigate the flood risks entirely.

Regulatory Framework
The regulatory framework in place during 2017 was the 2014 Plan, implemented by the International Joint Commission (IJC), which is responsible for overseeing water levels in the Great Lakes and St. Lawrence River. The 2014 Plan was designed to more naturally mimic the highs and lows of water levels to benefit aquatic life and wetlands, reduce the frequency of extreme water levels, and provide more reliable hydroelectric power production, while still considering the interests of shipping, recreational boating, and shoreline communities.

Regulated Outflows
In response to the extreme conditions of 2017, the International Lake Ontario-St. Lawrence River Board, which operates under the guidance of the IJC, adjusted outflows from Lake Ontario through the Moses-Saunders Power Dam, located on the St. Lawrence River between Cornwall, Ontario, and Massena, New York. The adjustments sought to balance the need to alleviate flooding in upstream areas around Lake Ontario with the risk of exacerbating downstream flooding in the Montreal area and other parts of the St. Lawrence River.

Despite these efforts, the extraordinary volume of water made it impossible to prevent flooding entirely. The board maximized outflows when conditions allowed; however, limitations were imposed by the need to maintain safe navigation conditions in the St. Lawrence Seaway and to minimize flooding impacts in downstream communities.

Impact and Response
The 2017 high water event led to widespread flooding and damage along the shores of Lake Ontario and the St. Lawrence River, prompting discussions about the effectiveness of the 2014 Plan and the need for adaptive management strategies in the face of climatic variability and extreme weather events. It underscored the challenges of managing a complex, binational water system in a changing climate, emphasizing the need for continued collaboration between the United States and Canada, as well as ongoing research and monitoring to improve future water management decisions.

\bibitem{AI3}
OpenAI ChatGPT 4\\
Query1: Give me some details about the specific regulations IJC make to control the water.\\
Output: The International Joint Commission (IJC), through Plan 2014, introduced a set of regulations aimed at managing the water levels and flows in Lake Ontario and the St. Lawrence River more closely to natural variations, while also considering human and environmental needs. Here are some specific aspects of how the IJC seeks to control water through these regulations:

1. Flow Targets and Limits

Plan 2014 establishes flow targets for the Moses-Saunders Power Dam, located on the St. Lawrence River, which is the primary control structure for outflows from Lake Ontario. These targets are designed to mimic natural water level variations more closely than the previous regulation plan, within the constraints of weather conditions and water supply.

2. Trigger Levels

The plan introduces "trigger levels" for water levels on Lake Ontario. When water levels reach these predefined trigger points, they activate certain adjustments in the outflows to either mitigate flooding when levels are too high or increase water levels when they are too low, especially during drought conditions. This approach aims to balance the needs of coastal communities, commercial navigation, hydropower generation, and environmental considerations.

3. Criteria for Deviations

Plan 2014 provides criteria for deviations from the standard flow regulations under specific conditions, such as extreme weather events or significant hydrological changes. These criteria allow for flexibility in response to unforeseen or exceptional circumstances, ensuring that the management of water levels can adapt to protect against extreme high or low water levels that could harm the environment or human activities.

4. Environmental Conditions

A significant part of the regulation focuses on improving environmental conditions along the shoreline and in wetlands. By allowing more natural fluctuations in water levels, the plan supports the restoration of habitats for fish, birds, and other wildlife. This aspect of the plan reflects a shift towards recognizing and valuing the ecological health of the Great Lakes-St. Lawrence River system.

5. Monitoring and Review

The plan includes provisions for ongoing monitoring and review of its impacts on water levels, environmental health, and economic factors. This ensures that the plan remains effective and can be adjusted as necessary based on actual outcomes, scientific understanding, and changing climate conditions.

Plan 2014 represents a complex regulatory framework that attempts to balance diverse interests and environmental sustainability in the management of Lake Ontario and the St. Lawrence River. Its implementation reflects an adaptive management approach, recognizing the challenges of climate variability and the need for flexibility in water resources management.

Query 2: Can you give me some website to know more about the detail of Plan 2014 and its future influence?

Output: Based on a quick search, here's what I found.

To learn more about Plan 2014 and its future influence on the Great Lakes and surrounding regions, the following websites offer detailed information and resources:

International Joint Commission (IJC) - Plan 2014 
Overview: The IJC's official website provides comprehensive information about Plan 2014, including its objectives, implementation details, and updates on its impact and future reviews. It's a primary source for official documents, reports, and latest news related to Plan 2014.

Plan 2014 | International Joint Commission
Great Lakes-St. Lawrence River Adaptive Management (GLAM) Committee: The GLAM Committee conducts reviews and assessments of Plan 2014, offering detailed reports on its effectiveness, environmental impacts, and adaptive management strategies in response to changing conditions.

Expedited Review of Plan 2014
Great Lakes Now: A platform that provides news and information about the Great Lakes, including coverage on Plan 2014. It offers insights into the plan's reception, its effects on local communities, and ongoing discussions about water management and environmental health in the Great Lakes region.

Great Lakes Now
These resources offer a wealth of information for understanding the complexities of Plan 2014, its regulatory mechanisms, environmental impacts, and the ongoing debates surrounding its future direction.

\bibitem{AI4}
OpenAI ChatGPT Classic\\
Query: make this more brief: 2017 was marked by exceptionally high precipitation in the Lake Ontario-St. Lawrence River basin during the spring months, contributing to record-breaking water levels in Lake Ontario. This was compounded by the runoff from a substantial snowpack in the upper Great Lakes and the Ottawa River basin, which flows into the St. Lawrence River upstream of Montreal. These natural conditions exceeded the capacity of any regulatory plan to mitigate the flood risks entirely.

The regulatory framework in place during 2017 was the 2014 Plan, implemented by the International Joint Commission (IJC), which is responsible for overseeing water levels in the Great Lakes and St. Lawrence River. The 2014 Plan was designed to more naturally mimic the highs and lows of water levels to benefit aquatic life and wetlands.

The 2017 high water event led to widespread flooding and damage along the shores of Lake Ontario and the St. Lawrence River, prompting discussions about the effectiveness of the 2014 Plan and the need for adaptive management strategies in the face of climatic variability and extreme weather events. It underscored the challenges of managing a complex, binational water system in a changing climate, emphasizing the need for continued collaboration between the United States and Canada, as well as ongoing research and monitoring to improve future water management decisions.\\
Output: In 2017, the Lake Ontario-St. Lawrence River basin experienced unprecedented spring precipitation and snowpack runoff, leading to historic high water levels and overwhelming regulatory flood mitigation efforts. The governing 2014 Plan by the International Joint Commission aimed to balance water levels for environmental and human use but faced challenges under these extreme conditions. The resulting severe flooding triggered a reevaluation of this plan, highlighting the difficulties of managing such a vast, cross-border water system amid climate change. It stressed the importance of adaptive management, international cooperation, and enhanced research for future resilience.

\end{ReportAiUse}

\end{document}